\documentclass[11pt]{article}
\usepackage{amscd}
\usepackage{verbatim}
\usepackage{amssymb}
\usepackage{amsmath}
\usepackage{amsthm}
\newcommand{\id}{\mathord{{\mathrm 1}\kern-0.27em{\mathrm I}}\kern0.35em}

\newcommand{\Rbb}{\mathbb{R}}

\begin{document}
\vskip 0 true cm \flushbottom
\begin{center}
\vspace{24pt} { \large \bf On long-time existence for the flow
of static metrics with rotational symmetry} \\
\vspace{30pt}
{\bf L Gulcev},$^{\dag}$ \footnote{lgulcev@math.ualberta.ca}
{\bf TA Oliynyk},$^{\ddag}$
\footnote{todd.oliynyk@sci.monash.edu.au} and
{\bf E Woolgar}$^{\dag}$ \footnote{ewoolgar@math.ualberta.ca} %%

\vspace{24pt} %%
{\footnotesize $^\dag$ Dept of Mathematical and Statistical
Sciences,
University of Alberta,\\
Edmonton, AB, Canada T6G 2G1.\\
$^\ddag$ School of Mathematical Sciences, Monash University,
Clayton, Vic 3800, Australia}
\end{center}
\date{\today}
\bigskip

\begin{center}
{\bf Abstract}
\end{center}
\noindent B List has proposed a geometric flow whose fixed points
correspond to solutions of the static Einstein equations of general
relativity. This flow is now known to be a certain Hamilton-DeTurck
flow (the pullback of a Ricci flow by an evolving diffeomorphism) on
${\mathbb R}\times M^n$. We study the ${\rm SO}(n)$ rotationally
symmetric case of List's flow under conditions of asymptotic
flatness. We are led to this problem from considerations related to
Bartnik's quasi-local mass definition and, as well, as a special
case of the coupled Ricci-harmonic map flow. The problem also occurs
as a Ricci flow with broken ${\rm SO}(n+1)$ symmetry, and has arisen
in a numerical study of Ricci flow for black hole thermodynamics.
When the initial data admits no minimal hypersphere, we find the
flow is immortal when a single regularity condition holds for the
scalar field of List's flow at the origin. This regularity condition
can be shown to hold at least for $n=2$. Otherwise, near a
singularity, the flow will admit rescalings which converge to an
${\rm SO}(n)$-symmetric ancient Ricci flow on ${\mathbb R}^n$.

%\maketitle

\setcounter{equation}{0}
\newpage

\section{Introduction}
\setcounter{equation}{0}

\subsection{List's flow}

\noindent Many of the most exciting recent developments in geometric
analysis have arisen from the study of geometric flow equations.
Among the most prominent examples, the Ricci flow has yielded a
proof of the Poincar\'e and Thurston conjectures \cite{Perelman1,
Perelman, CZ, MT} and the diffeomorphic $\frac14$-pinched sphere
theorem \cite{BS}, while the inverse mean curvature flow has been
used to prove the Riemannian Penrose conjecture \cite{HI1, HI2}. The
latter has important consequences in physics, prompting the question
of whether other geometric flow problems might also arise from
physics.

In physics one is often led to consider a metric of Lorentzian
signature. Many geometric flow equations in Riemannian geometry are
second-order parabolic (or at least quasi-parabolic), and therefore
they can be studied with the powerful tools of the maximum principle
and entropies. The problem is that this power is generally lost in
passing to Lorentzian signature.

However, Riemannian metrics can arise in physics problems, as the
case of the Riemannian Penrose conjecture shows. To see how this can
occur, consider the class of {\it stationary spacetimes} in general
relativity. These spacetimes have metrics that admit a timelike
Killing vector field. Examples are the G\"odel and exterior Kerr
metrics \cite{HE}, in which there are preferred observers who view
the time evolution of spacetime as nothing more than constant
rotation. Quotienting a stationary spacetime by the isometry
generated by the timelike Killing vector field, one obtains a
Riemannian metric on the base manifold, which is smooth if the
spacetime has no closed timelike curves and if the Killing orbits
are complete \cite{Harris}.

If the timelike Killing field is hypersurface-orthogonal, we arrive
at the class of {\it static spacetimes}, among which are the
exterior Schwarzschild and flat Minkowski spacetimes. For static
spacetimes, the rotation vanishes and the spacetime metric splits as
a warped product of a one-dimensional fibre over a Riemannian base
manifold. In particular, a static spacetime metric on ${\mathbb
R}\times M^n$, $n>2$, can be written as
\begin{equation}
ds^2=g_{\mu\nu}dx^{\mu}dx^{\nu}=-e^{2u}dx^0 dx^0
+e^{-\frac{2u}{n-2}}g_{ij}dx^idx^j,
\quad \frac{\partial u}{\partial x^0}=0,
\quad \frac{\partial g_{ij}}{\partial x^0}=0.
\label{eq1.1}
\end{equation}
Here $g_{\mu\nu}$ is a metric on ${\mathbb R}\times M^n$, while
$g_{ij}$ is a metric on $M^n$. Coordinates on $M^{n+1}$ are
$(x^{\mu})=(x^0,x^i)$, so Greek indices run over one more value
than Roman ones.

The vacuum Einstein equation is the condition that the metric
(\ref{eq1.1}) has vanishing Ricci curvature. If we apply the vacuum
Einstein equation to the metric $g_{\mu\nu}$ in (\ref{eq1.1}), we
obtain the equations
\begin{eqnarray}
&&R_{ij}-\left ( \frac{n-1}{n-2}\right )\nabla_i u \nabla_j u = 0
\ , \label{eq1.2}\\
&&\Delta u = 0 \ , \label{eq1.3}
\end{eqnarray}
for the Ricci tensor of the base $(M^n,g_{ij})$. Here $\Delta u :=
g^{ij}\nabla_i\nabla_j u$ is the Laplacian constructed from the
connection $\nabla_i$ compatible with $g_{ij}$ (note our signature
choice for the Laplacian). We note that equation (\ref{eq1.3}) is
redundant in that it can be derived from (\ref{eq1.2}) using the
contracted second Bianchi identity. It is therefore merely an
integrability condition for (\ref{eq1.2}).

%DEFINITION 1.1%%%%%%%%%%%%%%%%%%%%%%%%%%%%%%%%%%%%%%%%%%
\bigskip
\noindent {\bf Definition 1.1.} Equations (\ref{eq1.2}, \ref{eq1.3})
(or (\ref{eq1.2}) alone) are known as the {\it static vacuum
Einstein equations}. Solutions $g_{ij}$ are called {\it static
vacuum metrics}.
\bigskip

B List \cite{List1, List2}, in his PhD thesis under the direction of G
Huisken, presented a system of flow equations whose fixed points
solve the static vacuum Einstein equations.

%DEFINITION 1.2%%%%%%%%%%%%%%%%%%%%%%%%%%%%%%%%%%%%%%%%%%
\bigskip
\noindent {\bf Definition 1.2.} {\it List's system of equations} is
the system
\begin{eqnarray}
\frac{\partial g_{ij}}{\partial t} &=& -2 \left ( R_{ij}
-k_n^2 \nabla_i u \nabla_j u \right )\ , \label{eq1.4}\\
\frac{\partial u}{\partial t} &=& \Delta u\ . \label{eq1.5}
\end{eqnarray}
\bigskip

\noindent Note that $t$ is the flow parameter, not the time
coordinate in the spacetime metric (which we denote by $x^0$). The
metrics $g_{ij}(t;x)$ are a family of Riemannian metrics on an
$n$-manifold, $u(t;x)$ are a family of functions, and $k_n$ is an
arbitrary constant. When $k_n = \sqrt{\frac{n-1}{n-2}}$, the fixed
points of List's system are exactly the static vacuum metrics,
together with a harmonic function $u$. However, we will keep $k_n$
as an arbitrary constant (which obviously can be absorbed in $u$) so
that we can consider all $n\ge 2$.

A particularly useful equation easily derived from (\ref{eq1.4},
\ref{eq1.5}) is
\begin{equation}
\frac{\partial}{\partial t} \vert \nabla u \vert^2
= \Delta \vert \nabla u \vert^2-2\vert \nabla \nabla u \vert^2
-2k_n^2 \left ( \vert \nabla u \vert^2 \right )^2\ . \label{eq1.6}
\end{equation}

It is now realized that List's system of flow equations is in fact a
certain Hamilton-DeTurck flow in one higher dimension; i.e., List's
system is really a Ricci flow, modified by pulling back along an
evolving diffeomorphism (e.g., \cite{Li}). This does not make List's
flow any less interesting however. Recall that Hamilton-DeTurck flow
is given by
\begin{equation}
\frac{\partial g_{\mu\nu}}{\partial \lambda} = -2R_{\mu\nu} +
\pounds_X g_{\mu\nu}\ , \label{eq1.7}
\end{equation}
where $X$ is a vector field. To obtain List's system, choose
\begin{eqnarray}
g_{\mu\nu}dx^{\mu}dx^{\nu}&=&e^{2k_n u}d\tau^2+g_{ij}dx^idx^j\ ,
\label{eq1.8} \\
X&=&-\left (g^{ij}\nabla_i u\right )\frac{\partial}{\partial x^j}
\ . \label{eq1.9}
\end{eqnarray}
Note that the $g_{\mu\nu}$ in (\ref{eq1.8}) differs from that in
(\ref{eq1.1}).

It is the purpose of this paper to study the long-time existence
properties of solutions of List's system of equations that evolve
from a complete, asymptotically flat and rotationally symmetric
initial pair $(g(0),u(0))$, subject to the restriction that $g(0)$
does not admit a minimal hypersphere.

It will be convenient to choose a certain coordinate system
throughout the flow which will enable us to exploit the initial
absence of minimal hyperspheres. This will require that we work with
a DeTurck version of List's equations; i.e., that we pull back by a
further diffeomorphism on the base manifold. The DeTurck version of
equations (\ref{eq1.4}--\ref{eq1.6}) is
\begin{eqnarray}
\frac{\partial g_{ij}}{\partial t} &=& -2 \left ( R_{ij} -k_n^2
\nabla_i u \nabla_j u \right ) +\pounds_X g_{ij}\ , \label{eq1.10}\\
\frac{\partial u}{\partial t}& = & \Delta u + \pounds_X u
\, \label{eq1.11}\\
\frac{\partial}{\partial t} \left ( \vert \nabla u \vert^2\right )
&=& \Delta \vert \nabla u \vert^2-2\vert \nabla \nabla u \vert^2
-2k_n^2 \left ( \vert \nabla u \vert^2 \right )^2 +\pounds_X \left (
\vert \nabla u \vert^2\right ) \ , \label{eq1.12}
\end{eqnarray}
where the vector field $X$ generates the aforementioned
diffeomorphism, and it is this system that we will work with
directly.

\subsection{Motivations}

\noindent List's flow appears as a relatively simple case of the
Ricci-harmonic map flow, which has been studied in \cite{RM}. The
general form of this flow is
\begin{eqnarray}
\frac{\partial g_{ij}}{\partial t}
&=& -2R_{ij}+G_{ab}\nabla_i u^a \nabla_j u^b\ , \label{eq1.13}\\
\frac{\partial u^a}{\partial t}&=&\Delta u^a +g_{ij}\Gamma^a_{bc}
\nabla^i u^b \nabla^j u^c\ , \label{eq1.14}
\end{eqnarray}
where the $u^a$ are embedding functions $u^a:(M^n,g_{ij})
\hookrightarrow ({\cal M}^m,G_{ab})$ mapping one Riemannian manifold
to another and the $\Gamma^a_{bc}$ are the coefficients of the
Levi-Civit\`a connection of $G_{ab}$. In the case that ${\cal M}^m
={\mathbb R}$, $u^a=u$, and $G_{ab}=2k_n^2$, these equations reduce
to List's flow. Thus, List's flow is the special case of the coupled
Ricci-harmonic map flow where the target space $({\cal M}^m,G_{ab})$
is the real line.

List's flow with rotational symmetry also appears if rotational
symmetry is broken in Ricci flow in one more dimension. Rotationally
symmetric, asymptotically flat Ricci flow in dimension $n\ge 3$ was
studied in connection with a conjecture in string theory regarding
the limiting behaviour of ADM mass as the flow converges \cite{GHMS}
(see \cite{Ivey} for an earlier study and see \cite{Wu} for the
$n=2$ case). This is now well-understood, at least in the absence of
minimal hyperspheres \cite{OW} (see also \cite{DM}). It is
interesting to ask how this understanding is modified if the
rotational symmetry is broken down to a subgroup. By the above
correspondence between flows, we see that List's flow with
rotational symmetry can be thought of as Ricci flow on an
$(n+1)$-manifold with ${\mathbb R}\times{\rm SO}(n)$ symmetry, which
has $\frac{n(n-1)}{2}+1$ generators, or $n-1$ fewer generators than
the $\frac{n(n+1)}{2}$ generators of full rotational symmetry in
$(n+1)$-dimensions.

However, there is a another reason to study this system, which may
prove to be the most compelling. List's equations were conceived as
a tool to address conjectures about static metrics in general
relativity \cite{Huisken}. We briefly discuss two of these.

We recall Bartnik's quasi-local definition of mass \cite{Bartnik}.
In an $(n+1)$-dimensional spacetime, consider a moment of time
symmetry (i.e., a spacelike hypersurface with zero extrinsic
curvature) and in it a bounded $n$-dimensional region $B$. Consider
all asymptotically flat Riemannian $n$-manifolds $N$ of non-negative
scalar curvature $R\ge 0$ into which $B$ can be isometrically
embedded (smoothly in the interior of $B$), such that the induced
metric and mean curvature must match from both sides of $\partial
B$. Further assume that $N$ has no stable minimal sphere lying
outside the image of $B$. Then $N$ is called an {\it admissible
extension} of $B$. By the positive mass theorem $N$ has nonnegative
Arnowitt-Deser-Misner (ADM) mass. Consider all such admissible
extensions of $B$ and take the infimum of the ADM masses. This
infimum is the {\it Bartnik mass} $m_B$ of the region $B$. It is
clearly nonnegative.

What is not so clear from the definition is whether the mass ever
differs from zero. This led Bartnik to make the following conjecture
which, if true, would guarantee that the Bartnik mass is nontrivial:

%Static Minimization Conjecture %%%%%%%%%%%%%%%%%%%%%%%%%%%%%%%%%%%
\bigskip
\noindent{\bf Static Minimization Conjecture (Bartnik):} The infimum
is a minimum, and is realized as the ADM mass of a solution of the
static vacuum Einstein equations.
\bigskip

\noindent Huisken and Ilmanen \cite{HI2} have since shown by other
methods that Bartnik's mass is nonzero except when $B$ is a domain
in flat space, thus proving the nontriviality of the Bartnik mass.
However the static minimization conjecture has remained open up to
now.\footnote
{However, as we were preparing the final draft of this manuscript, a
preprint appeared \cite{AK} announcing a proof that for any bounded
3-dimensional spatial region whose boundary has positive Gauss
curvature, there exists an extension satisfying the static Einstein
equations with suitable boundary conditions (Bartnik's {\it
geometric conditions}).}

One strategy to address this conjecture may be to choose one
admissible, asymptotically flat extension of $B$ and use it as the
initial condition for List's flow. Boundary conditions, such as
Bartnik's {\it geometric conditions} (\cite{Bartnik2}) that fix the
boundary mean curvature and induced metric, must also be imposed at
$\partial B$. The idea is then to use the flow to produce a
mass-minimizing sequence which converges to a fixed point, hence to
a static metric.

A test case would be to employ this strategy on ${\mathbb R}^n$,
with no inner boundary. Ideally this would produce sequences of
metrics that converge to flat space.\footnote
{Note that the manner in which List's flow would produce
mass-minimizing sequences will be similar to that of the Ricci flow.
There the mass remains constant throughout the flow but will jump to
a minimizing value in the limit as $t\to\infty$, while various
measures of the quasi-local mass within bounded regions flow
smoothly toward minimizing values \cite{OW, GHMS, DM}.}
Huisken and Ilmanen \cite{HI2} have discussed such mass minimizing
sequences, and suggest a more complicated view.

\bigskip
\noindent{\bf Conjecture (Huisken and Ilmanen).} Suppose $(M,g_i)$
is a sequence of asymptotically flat, mass-minimizing, non-negative
scalar curvature 3-metrics tending to zero mass. Then there is a set
$Z_i$ such that $|\partial Z_i|\to 0$ and $(M\backslash Z_i,g_i)$
has a flat Gromov-Hausdorff limit.\footnote
{at least, if the scalar curvature of each $g_i$ is nonnegative.}

\bigskip
\noindent This foresees that an obstruction to convergence may arise
in rotationally symmetric, asymptotically flat List flow, in the
form of a locally collapsed long, thin tube growing at the
origin.\footnote
{We expect that if collapse occurred elsewhere, the rotational
symmetry would force this to be preceded by formation of a minimal
surface. But we will show that the absence of an initial minimal
surface implies that none can form later.}

There is some numerical evidence in favour of convergence to flat
space. As part of a study motivated by black hole thermodynamics,
Headrick and Wiseman \cite{HW} examined Ricci flow manifolds-with
boundary with ${\rm U}(1)\times {\rm SO}(3)$ symmetry, including
$S^1\times{\mathbb B}^3$ where ${\mathbb B}^3$ denotes a 3-ball in
${\mathbb R}^3$. They thus had a finitely distant spatial boundary
and imposed a Dirichlet condition there. On $S^1\times{\mathbb B}^3$
they found convergence to flat $S^1\times {\mathbb R}^3$, known in
the physics literature as ``hot flat space''.

We therefore undertook a study of the long-time existence properties of
this flow, with a view to shedding analytical light on these
conjectures and numerical results.

\subsection{Overview and main results}

\noindent Even with the restriction to rotational symmetry, the
above conjectures are not easy to address, and our results are only a
starting point. We prove the following:

\bigskip
\noindent{\bf Theorem 1.3.} {\sl Let $\left ( {\tilde
g}_{ij}(r),{\tilde u}(r)\right )$ be asymptotically flat and
rotationally symmetric initial data for the system of flow equations
(\ref{eq1.10}--\ref{eq1.12}) on ${\mathbb R}^n$ such that the metric
${\tilde g}_{ij}(r)$ admits no minimal hypersphere. Then this system
of equations has an asymptotically flat, rotationally symmetric
solution $\left ( g_{ij}(t,r), u(t,r) \right )$ on $[0,T_M)\times
[0,\infty)$ for some maximal time of existence $T_M\in (0, \infty ]$.
No minimal hypersphere forms at any $t<\infty$. Furthermore,
\begin{enumerate}
\item [{\it (i)}] If $n=2$, then the flow is immortal
    ($T_M=\infty$).
\item [{\it (ii)}] If $n\ge 3$, and if there is a function
    $F:[0,\infty)\to (0,\infty)$ such that  $\frac{1}{r}|\nabla
    u|_{(t,r)}\le F(t)$, then the flow is immortal.
\end{enumerate}
}
\bigskip

In the case where the flow fails to exist, we can go some short
distance towards analyzing the kind of singularity that develops.
List has shown in his thesis \cite{List1} that where the flow fails
to exist, the norm of the Riemann tensor diverges. We therefore
borrow the following definition from the Ricci flow.

\bigskip
\noindent{\bf Definition 1.4.} For $T_M$ the maximal time of existence
of the Ricci flow, an {\it essential blow-up sequence} $(t_k,x_k)$
is a sequence of spacetime points such that $t_k \nearrow T_M$ and
$\sup_{[0,t_k]\times [0,\infty)}|{\rm Riem}|(t,r)\le C |{\rm
Riem}(t_k,r_k)|=:B_k$ for some constant $C\ge 1$.

\bigskip
\noindent {\bf Theorem 1.5.} {\sl Let $(g(t),u(t))$ be a
rotationally symmetric solution of (\ref{eq1.10}--\ref{eq1.12})
developing from initial data as in Theorem 1.3, with maximal time of
existence $T_M<\infty$, and let $({\bar g}(t),{\bar u}(t))
=\varphi_t^*(g,u)$ be the corresponding solution of (\ref{eq1.4},
\ref{eq1.5}). Let $(t_k,r_k)$, $t_k\nearrow T_M$, be an essential
blow-up sequence for $\left ({\bar g}(t),{\bar u}(t)\right )$. Set
$B_k:=\overline{ |{\rm Riem}|}(t_k,r_k)$ and define the rescalings
\begin{eqnarray}
g_{(k)(s)}&:=&B_k\cdot {\bar g}(t_k + s/B_k)\ , \nonumber\\
u_{(k)}(s)&:=&{\bar u}(t_k+ s/B_k)\ , \label{eq1.15}\\
s&\in& [-B_k(1+t_k),0]\ .\nonumber
\end{eqnarray}
Then there is a subsequence of the pointed sequence $({\mathbb R}^n,
g_{(k)}(s), u_{(k)}(s), r_k)$ which converges smoothly on all
compact subsets of $(-\infty,0]\times {\mathbb R}^n$ to $({\mathbb
R}^n, g_{(\infty)}(s), const , r_{\infty})$ with $({\mathbb R}^n,
g_{(\infty)}(s))$ a complete, ancient solution of the Ricci flow.
Sectional curvature in planes tangent to the orbits of the
rotational symmetry group is nonnegative, and so is the scalar
curvature. In $n=3$ dimensions, sectional curvature in planes
containing the radial vector is also nonnegative.}
\bigskip

Theorem 1.5 does not confirm the conjecture above of
Huisken-Ilmanen, as it leaves open the possibility that the flow is
non-collapsed below some finite scale at the singularity time. When
that occurs, then in any dimension, after rescaling, the resulting
limit would be noncollapsed below any scale and, in the $n=3$ case,
would have nonnegative sectional curvatures. By rotational symmetry,
the limit would then be a Bryant soliton for $n=3$.

In section 2 we discuss the notion of asymptotic flatness that we
use, and survey results of List on local existence and continuation,
making minor modifications where necessary. Section 3 discusses
rotational symmetry and its implications. It is in this section that
we state the evolution equations in the form that we use and define
the basic quantities whose flow we analyse in subsequent sections.
Section 4 contains estimates that are valid in arbitrary dimension
with no further assumptions beyond rotational symmetry and
asymptotic flatness. In section 5, we assume either that
$\frac{1}{r}|\nabla u|$ is bounded on any closed time interval or
that the dimension is $n=2$ (in which case it is shown in section 4
that $\frac{1}{r}|\nabla u|$ is bounded on closed time intervals).
Under either of these assumptions, we are then able to obtain all
further estimates required to show boundedness of sectional
curvatures on finite time intervals. The proofs of Theorems 1.3 and
1.5 then follow easily form these results. These proofs are given in
section 6.

Our sign and index conventions are as follows. As previously stated,
we take the (rough or scalar) Laplacian to be
$\Delta:=g^{ab}\nabla_a \nabla_b$. We define the curvature
$R^a{}_{bcd}x^by^cz^d:=\nabla_y\nabla_z x-\nabla_z\nabla_y x
-\nabla_{[y,z]}x$, and we write $R_{abcd}:=g_{ae}R^e{}_{bcd}$. The
Ricci tensor is $R_{bd}:=R^a{}_{bad}$. We endeavor where possible to
denote constants that bound a function $h$ (say) by $C^+_h$ for an
upper bound (i.e., to indicate that $h(t,r)\le C^+_h$ for all $t$)
and $C^-_h$ for a lower bound, though we sometimes deviate from this
practice for reasons of convenience.

%SECTION 2%%%%%%%%%%%%%%%%%%%%%%%%%%%%%%%%%%%%%%%%%%%%
\section{Preliminaries: Asymptotic flatness}
\setcounter{equation}{0}

%

%Asymptotic flatness%%%%%%%%%%%%%%%%%%%%%%%%%%%%%%%%%%%%%%%%
%\subsection{Asymptotic flatness}

\noindent The definition of asymptotic flatness, more properly
called local asymptotic flatness when $n=2$, can be formulated on
any Riemannian manifold with dimension $n\ge 2$ which admits the
notion of an asymptotic end. However, since we work on ${\mathbb
R}^n$, complete generality is not necessary here, though it can be
achieved with minor changes to the formulation below. On the other
hand, our results will hold with a much more general notion of
asymptotic flatness for $n\ge 3$ than the usual notion.

To begin, we let
\begin{equation}
e_{ij} = \begin{cases} \delta_{ij} & \text{for $n\geq 3$} \label{eq2.1}\\
  \delta_{ij}+a\frac{x_i x_j}{r^2} & \text{for $n=2$}
\end{cases}
\end{equation}
where $r = \sqrt{\sum_{i=1}^n (x^i)^2}$ and, in dimension $n=2$,
$x_i=\delta_{ij}x^j$, $a>-1$, and the deficit angle of the flat cone
metric is $2\pi\left (1-\frac{1}{\sqrt{1+a}}\right )$.

%%EW: The following replaces the definition of asymptotic flatness.

Following \cite{List1, List2}, we define

\bigskip
\noindent{\bf Definition 2.1.} For $n\ge 3$, $(M,g,u)$ is {\it
asymptotically flat (of order one)} if there is a compact subset
$K\subset M$ such that $M_K:=M\backslash K$ is diffeomorphic to
${\mathbb R}^n\backslash B_1(0)$ where $B_1(0)$ is the Euclidean
unit ball and, on $M_K$, $(g,u)$ satisfies
\begin{eqnarray}
\left \vert g_{ij}-e_{ij}\right \vert
&\le& C_0/r \ , \label{eq2.2} \\
\left \vert \partial_k g  \right \vert &\le& C_k/r^{k+1}\ ,
\ k=1,2,3, \label{eq2.3}\\
\left \vert u\right \vert &\le& D_0/r \ , \label{eq2.4} \\
\left \vert \partial_k u \right \vert &\le& D_k/r^{k+1}\ ,
\ k=1,2,3, \label{eq2.5}
\end{eqnarray}
where $C_k$, $D_k$ are constants ($k\in \{ 0,1,2,3\}$),
$r^2=x_1^2+\dots +x_n^2$ with the $x_i$ being Cartesian coordinates
for $e_{ij}$, and $\partial_k$ is the Cartesian coordinate
derivative.
\bigskip

We choose to work in this class for three reasons. The first is that
a local-in-time existence theorem within this class is already
available \cite{List1, List2}. The second is that these fall-off
conditions are well-suited to the arguments in subsequent sections.
The third is that our results are, in fact, not very sensitive to
the precise choice of definition of asymptotic flatness. We
therefore settle on a convenient choice rather than the most general
one, for which the preliminaries would be a greater
distraction.\footnote
{For example, a definition based on weighted Sobolev spaces was used
for a similar problem in \cite{OW}. One can augment that definition
by including a condition that $u$ lie in a weighted Sobolev space
$H_{\delta}^{k}$ with $k$, the number of derivatives in the Sobolev
norm, chosen such that $k>3+n/2$ and $\delta$, the exponent in the
weight factor $r^{\delta}$, any negative number. Then local-in-time
existence can be obtained in this class, by modifying the argument
in \cite{OW} and papers cited therein.}

Having said that, we note that when $n>3$ this definition is in fact
much {\it weaker} than most, which tend to require the difference
between the metric and $e_{ij}$ to decay at ${\cal O}(1/r^{n-2})$.
However, for $n=3$, the present definition is stronger than
necessary.

\bigskip
\noindent{\bf Proposition 2.2 (List).} {\sl Let $({\hat g},{\hat
u})$ be asymptotically flat. Then there exists a $T>0$ such that
$(g(t),u(t))$ solves (\ref{eq1.4}, \ref{eq1.5}) for all $0\le t <T$
and $g(0)={\hat g}$, $u(0)={\hat u}$. Furthermore, $(g(t),u(t))$ is
asymptotically flat for all $0\le t <T$.}

\bigskip
\noindent{\bf Proof.} See \cite{List1}, Theorems 3.12 and 9.7, or
\cite{List2}, Theorems 4.1 and 8.6. \qed

\bigskip
\noindent{\bf Remark 2.3.} In fact, List gives a detailed proof of
Theorem 2.1 assuming stronger asymptotic flatness conditions and
then notes that the proof obviously goes through as well for
asymptotic flatness conditions which agree with those above when
$n\ge 3$. This is clearly the case, and furthermore it is also the
case for $n=2$, with $e_{ij}$ used in place of $\delta_{ij}$ at one
step in the proof (equation (9.7) of \cite{List1} or equation (8.5)
of \cite{List2}).
\bigskip

It is now possible to state a criterion for the flow to exist for
all future time, in the form of a continuation principle which
states that, as with Ricci flow, the flow can be continued beyond
$t=T$ unless the norm of the Riemann curvature diverges there.

\bigskip
\noindent {\bf Proposition 2.4.} {\sl Let $(\hat{g},{\hat u})$ be a
asymptotically flat initial data. Then the system (\ref{eq1.4},
\ref{eq1.5}), with the initial conditions $g(0)=\hat{g}$,
$u(0)=\hat{u}$ has a unique solution on a maximal time interval
$0\leq t < T_M \leq \infty$. If $T_M < \infty$ then
\begin{equation}
\limsup_{t\nearrow T_M}
%\underset{\;t\nearrow T_M} {\overline{\lim}\;\;}
\sup_{x\in\Rbb^n}|{\rm Rm}(t,x)|_{g(t,x)} = \infty \ . \label{eq2.6}
\end{equation}
Moreover, for any $T \in [0,T_M)$, if $K= \sup_{0\leq t\leq T}
\sup_{x\in\Rbb^n}|{\rm Rm}(t,x)|_{g(t,x)}$, and $C=\sup_{x\in
\Rbb^n}|\nabla \hat{u}(x)|_{g(x)}^2$, then}
\begin{equation} e^{-(2nK+4C)T}\hat{g} \leq g(t) \leq e^{(2nK+4C)T}\hat{g}
\quad \text{for all\ } t\in [0,T]\ . \label{eq2.7}
\end{equation}
\bigskip

\noindent{\bf Proof.} List gives a partial proof for complete
manifolds (Theorem 3.22 of \cite{List1}) and a full proof for closed
manifolds (Theorems 3.11 and 6.22 of \cite{List1}). The full proof
uses the closedness of the manifold only to invoke the maximum
principle for nonnegative scalar functions; in particular, norms of
$\nabla u$, ${\rm Riem}$, and derivatives thereof. By Proposition
2.2 and equations (\ref{eq2.2}--\ref{eq2.5}), each such quantity
tends to zero as $r\to\infty$, $0\le t\le T$, and thus the maximum
principle applies to these functions on complete manifolds as well,
provided the initial data obey asymptotic flatness. \qed

\bigskip

%SECTION 3%%%%%%%%%%%%%%%%%%%%%%%%%%%%%%%%%%%%%%%%%%%%
\section{Rotational Symmetry}
\setcounter{equation}{0}

\subsection{The coordinate system}
\noindent Now assume the flow that solves (\ref{eq1.4}, \ref{eq1.5})
evolves from rotationally symmetric $C^2$ initial data $({\hat
g},{\hat u})$ with
\begin{equation}
d{\hat s}^2= {\hat g}_{ij}dx^idx^j=a^2(\rho)d\rho^2+\rho^2
d\Omega^2\ ,\label{eq3.1}
\end{equation}
where $d\Omega^2$ is the constant curvature ${\rm sec}=1$ metric on
the $(n-1)$-sphere. We take $a(0)=1$, $a'(0)=0$.

Ricci flow preserves isometries. List's flow, in turn, preserves
symmetries of the pair $(g,u)$ (i.e., isometries of $g$ that commute
with $u$). Combining this fact with Proposition 2.4, then there will
be a maximal time of existence $T_M\in (0,\infty]$, a coordinate
system in which (\ref{eq2.1}--\ref{eq2.5}) hold, and coordinate
transformations taking the metric to a spherical coordinate system
$x^i=(\rho,\theta^A)$ (with $\theta^A$ the coordinates on the
$\rho=const$ spheres). In these coordinates, the flow is
\begin{eqnarray}
t &\mapsto& \left ( {\bar g}(t,\rho), {\bar u}(t,\rho)\right ) \nonumber \\
d{\bar s}^2 &=& {\bar g}_{ij}dx^idx^j=q^2(t,\rho)d\rho^2+h^2(t,\rho)
d\Omega^2\ .\label{eq3.2}
\end{eqnarray}
This metric solves (\ref{eq1.4}, \ref{eq1.5}). The coordinate
functions $q$ and $h$ are $C^2$ in $\rho$ and, in the one-sided
sense, $C^1$ in $t$ at $t=0$ and $\rho>0$, and are smooth in $t$ and
$r$ for all $r>0$ and $t\in (0,T_M)$.

Now introduce a new coordinate system at each time, obtained via
acting with the family of diffeomorphisms
\begin{equation}
\psi_t(\rho,\theta^A):=\left ( h(t,\rho), \theta^A \right )
=: (r,\theta^A) \ . \label{eq3.3}
\end{equation}
Note that since $h(0,\rho)=\rho$ then $\psi_0={\rm id}$, and also,
since $\frac{\partial h}{\partial \rho}(0,\rho)=1$, then for $T$
sufficiently small, there are (possibly $T$ dependent) constants
$C_{\partial h}^{\pm}(T)$ such that
\begin{equation}
0<C^-_{\partial h}(T) \le \frac{\partial h}{\partial \rho}(t,\rho)
\le C^+_{\partial h}(T) \label{eq3.4}
\end{equation}
whenever $0\le t \le T$ and $T<{\tilde T}$, where ${\tilde T}\le
T_M$ is defined to be the supremum of $T$-values for which
(\ref{eq3.4}) is true.

\bigskip
\noindent{\bf Proposition 3.1.} ${\tilde T}=T_M$.
\bigskip

\noindent{\bf Proof.} Given in subsection 4.2.2. \qed
\bigskip

\noindent We can now write the flow in ``area radius
gauge'' as
\begin{eqnarray}
t&\mapsto&\left ( g(t,r),u(t,r)\right ) \nonumber \\
g(t,r)&:=&\left ( \psi_t^{-1} \right )^* {\bar g}(t,\rho)
=f^2(t,r)dr^2+r^2d\Omega^2 \ , \label{eq3.5}\\
u(t,r)&:=&\left ( \psi_t^{-1} \right )^* {\bar u}(t,\rho)
=u\left ( t,h(t,\rho) \right ) \ , \nonumber
\end{eqnarray}
where
\begin{equation}
f(t,r):=\frac{q\left ( t, \rho(t,r) \right )}
{\frac{\partial h}{\partial \rho}
\left ( t, \rho(t,r) \right )}\ . \label{eq3.6}
\end{equation}

Let the generator of the family $\psi_t$ be written as
$X^j=g^{ij}\nabla_i\phi(t,r)$ for some scalar $\phi(t,r)$. Inserting
this and (\ref{eq3.5}) into (\ref{eq1.10}), we obtain the pair of
equations
\begin{equation}
\frac{\partial f}{\partial t}
=-\frac{(n-1)}{rf^2}\frac{\partial f}{\partial r}
+{k_n^2f\left |\nabla u \right \vert^2}
+\frac{1}{f}\frac{\partial^2\phi}{\partial r^2}
-\frac{1}{f^2} \frac{\partial f}{\partial r}
\frac{\partial \phi}{\partial r}\ , \label{eq3.7}
\end{equation}
\begin{equation}
0=\frac{r}{f^3}\frac{\partial f}{\partial r}
+(n-2)\left(1-\frac{1}{f^2}\right)
-\frac{r}{f^2}\frac{\partial \phi}{\partial r}
\ .\label{eq3.8}
\end{equation}
The latter yields
\begin{equation}
\frac{\partial \phi}{\partial r}
=\frac{1}{f}\frac{\partial f}{\partial r}
+\frac{(n-2)}{r}(f^2-1)\ , \label{eq3.9}
\end{equation}
We substitute this into (\ref{eq3.7}) and (\ref{eq1.11}) and define
\begin{equation}
z:=\frac{1}{f} \frac{\partial u}{\partial r}\label{eq3.10}
\end{equation}
so that $z^2=\vert \nabla u \vert^2$ in rotational symmetry. Then
equations (\ref{eq1.10}, \ref{eq1.12}) reduce to the system of
equations which we study herein, namely:

%DEFINITION 3.2%%%%%%%%%%%%%%%%%%%%%%%%%%%%%%%%%%%%%%%%%%
\bigskip
\noindent{\bf Definition 3.2.} The {\em rotationally symmetric flow
equations} are the system
\begin{eqnarray}
\frac{\partial f}{\partial t}&=&\frac{1}{f^2}
\frac{\partial^2 f}{\partial r^2}-\frac{2}{f^3}\left (
\frac{\partial f}{\partial r} \right)^2+\left (\frac{n-2}{r}
-\frac{1}{rf^2} \right)\frac{\partial f}{\partial r} \nonumber \\
&&-\frac{(n-2)}{r^2f}(f^2-1)+k^2_nfz^2\ ,
\label{eq3.11}\\
\frac{\partial z}{\partial t}&=&\frac{1}{f^2}
\frac{\partial^2 z}{\partial r^2}+\left[\frac{1}{rf^2}
+\frac{n-2}{r} \right]
\frac{\partial z}{\partial r}-\left[\frac{n-1}{r^2f^2}
+k^2_nz^2\right]z \label {eq3.12}\ .
\end{eqnarray}
\bigskip

Next, define
\begin{eqnarray}
\lambda_1(t,r)&:=&\frac{1}{rf^3}\frac{\partial f}{\partial r}
\ ,\label{eq3.13}\\
\lambda_2(t,r)&:=&\frac{1}{r^2}\left ( 1- \frac{1}{f^2} \right )
\ . \label{eq3.14}
\end{eqnarray}

\bigskip
\noindent{\bf Lemma 3.3.} {\sl When $n=2$, $\lambda_1$ is the Gauss
curvature. When $n\ge 3$, $\lambda_1$ is the sectional curvature in
planes containing $\frac{\partial}{\partial r}$ and $\lambda_2$ is
the sectional curvature in planes tangent to the $r=const$ spheres.
As well, we have}
\begin{eqnarray}
\left |{\rm Riem}\right \vert ^2 &=&R_{ijkl}R^{ijkl}=2(n-1)
\lambda_1^2+(n-1)(n-2)\lambda_2^2\ ,\label{eq3.15}\\
\frac{\partial \lambda_2}{\partial r}
&=& \frac{2}{r} \left ( \lambda_1 -\lambda_2 \right )\ . \label{eq3.16}
\end{eqnarray}

\bigskip

\noindent {\bf Proof.} The curvature interpretations of $\lambda_1$
and $\lambda_2$ follow from trivial computations, and then
(\ref{eq3.15}) follows immediately from rotational symmetry.
Equation (\ref{eq3.16}) is obvious (expand both sides) and is, in
fact, the second Bianchi identity. \qed

\bigskip
Note that (\ref{eq3.16}) shows that $\lambda_1=\lambda_2$ at the
origin and, more generally, at any spatial or spacetime extremum of
$\lambda_2$. Also note that, using Proposition 3.1. and Lemma 3.3,
we can adapt the continuation principle (Proposition 2.4) to the
area-radius gauge:

%Proposition 3.4%%%%%%%%%%%%%%%%%%%%%%%%%%%%%%%%%%%%%%%%%%
\bigskip
\noindent {\bf Proposition 3.4.}  {\sl If there exists a constant
$C_\lambda >0$ independent of $T_M$ such that
\begin{equation}
%\sup_{0<r<\infty}
|\lambda_1(t,r)| \leq C_\lambda \quad \text{ if $n=2$,} \label{eq3.17}
\end{equation}
or
\begin{equation}
%\sup_{0<r<\infty}\bigl(
|\lambda_1(t,r)|+|\lambda_2(t,r)|
%\bigr)
\leq C_\lambda \quad \text{ if $n\geq 3$,} \label{eq3.18}
\end{equation}
 for all $(t,r)\in [0,T_M)\times [0,\infty)$ then $T_M=\infty$.}

\bigskip

\noindent We shall eventually see that it suffices to bound
$\lambda_1$ or, equivalently, $R$ from above.

%SECTION 4%%%%%%%%%%%%%%%%%%%%%%%%%%%%%%%%%%%%%%%%%%%%
\section{A priori bounds}
\setcounter{equation}{0}

\noindent In this section and the next, we always assume that
$(g(t),u(t))$, $0\le t \le T$, is a solution of
(\ref{eq1.10}--\ref{eq1.12}) and that $(M,g(t),u(t))$ is
asymptotically flat (of order one; as in Definition 2.1) for all
$t\in [0,T]$. For now $T< {\tilde T}\le T_M$, but after we prove
Proposition 3.1 we will be able to set ${\tilde T}=T_M$.

\subsection{Elementary bounds on scalar quantities}

\noindent In his thesis \cite{List1, List2}, List shows that a
modified form of the usual Ricci flow lower bound on scalar
curvature of compact manifolds holds for List's flow. He also shows
that $|\nabla u|^2$ is bounded above by $const/(1+t)$ on a compact
manifold. These results are simple applications of the maximum
principle. In this section, we adapt the maximum principle to the
complete, asymptotically flat setting and obtain bounds on $R$ and
$|\nabla u|^2$ for asymptotically flat manifolds as corollaries.

%LEMMA 4.1%%%%%%%%%%%%%%%%%%%%%%%%%%%%%%%%%%%%%%%%%%%%
\bigskip
\noindent{\bf Lemma 4.1.} {\sl Let $\Psi$ be a solution of
\begin{equation}
\frac{\partial \Psi}{\partial t}\le \Delta \Psi
+ \nabla_Y \Psi - k^2 \Psi^2 \label{eq4.1}
\end{equation}
for some vector field $Y$ and constant $k$ on the domain
$D(T):=[0,T]\times {\mathbb R}^n \ni (t,x)$, such that $\Psi\to 0$
as $x\to\infty$. (i) If $\Psi(0,x)\le 0$ for all $x\in {\mathbb
R}^n$ then $\Psi(t,x)\le 0$ for all $(t,x)$ in $D(T)$, and otherwise
(ii) if $k\neq 0$ we have $\Psi \le C^+_{\Psi}/(1+t)$ for
$C^+_{\Psi}:=\max \left \{\frac{1}{k^2}, \sup_{x\in{\mathbb R}^n}
\Psi(0,x)\right \}$.}

\bigskip
\noindent{\bf Proof.} Consider the domain $D_{\epsilon}(T)
:=[0,T]\times B^n(1/\epsilon)$ for $B^n(1/\epsilon)$ the $n$-ball of
radius $1/\epsilon>0$ with respect to the Euclidean metric centred
at the origin in ${\mathbb R}^n$. To prove (i), we note that
(\ref{eq4.1}) implies that $\frac{\partial \Psi}{\partial t}\le
\Delta \Psi + \nabla_Y \Psi$, and so standard maximum principle
arguments show that $\Psi$ must take its maximum on the parabolic
boundary of $D_{\epsilon}$ (i.e., points where either $t=0$ or
$r=1/\epsilon$). Taking $\epsilon\to 0$, we see by the asymptotic
condition on $\Psi$ that either the supremum is $0$ or the supremum
occurs at $t=0$, in which case it is again zero by assumption.

To prove (ii), we define the function $G:=(1+t)\Psi$ and see that
from (\ref{eq4.1}) it obeys
\begin{equation}
\frac{\partial G}{\partial t}\le \Delta G + \nabla_Y G
+\frac{G}{1+t} \left [ 1-k^2 G \right ] \ .\label{eq4.2}
\end{equation}
If the maximum of $G$ on $D_{\epsilon}$ occurs at some $(t,x)$ in
the parabolic interior of $D_{\epsilon}$ (i.e., the complement of
the parabolic boundary), we see immediately from (\ref{eq4.2}) that
$G\le 1/k^2$, and so $\Psi\le \frac{1}{k^2(1+t)}$. If instead the
maximum of $G$ occurs on the boundary of $B^n(1/\epsilon)$ at some
$0<t<T$, then by taking $\epsilon$ small enough we see that the
maximum approaches zero. Alternatively, the maximum can occur at
$t=0$, and since $G(0,x)=\Psi(0,x)$ then in this case $\Psi(t,x)\le
\sup_x \frac{G(0,x)}{1+t}= \sup_x \frac{\Psi(0,x)}{1+t}$. \qed
\bigskip

This immediately yields Propositions 4.2 and 4.4.

%PROPOSITION 4.2%%%%%%%%%%%%%%%%%%%%%%%%%%%%%%%%%%%%%%%%
\bigskip
\noindent {\bf Proposition 4.2.} {\sl Let $(g,u)$ be an
asymptotically flat solution of the flow on the domain $[0,T]\times
{\mathbb R}^n$. $0<T<{\tilde T}$ Then for any $n\ge 2$
\begin{equation}
\left |\nabla u (t,r)\right \vert \le \frac{1}{\sqrt{1+t}}
C^+_{|\nabla u |}\ , \label{eq4.3}
\end{equation}
where $C^+_{|\nabla u |}$ is a constant depending only on the
initial data and $k_n$.}

\bigskip
\noindent{\bf Proof.} Equation (\ref{eq1.12}) yields
\begin{equation}
\frac{\partial}{\partial t} \left ( |\nabla u |^2 \right ) \le
\Delta \left ( |\nabla u |^2 \right )+\nabla_X \left (
|\nabla u |^2 \right ) -2k^2_n \left ( |\nabla u |^2 \right )^2
\ . \label{eq4.4}
\end{equation}
Apply Lemma 4.1 with $\Psi=|\nabla u |^2$ and $k^2=2k^2_n$. This
proves the proposition and shows that $C^+_{|\nabla u |}=\max \left
\{ \frac{1}{\sqrt{2k^2_n}}, \sup_x |\nabla u|_{(0,x)} \right
\}$.\qed

%REMARK 4.3%%%%%%%%%%%%%%%%%%%%%%%%%%%%%%%%%%%%%%%%%%%
\bigskip
\noindent{\bf Remark 4.3.} In rotational symmetry, we have $|z|=\left
\vert \frac1f \frac{\partial u}{\partial r}\right \vert =|\nabla u|$
and we write this bound as
\begin{equation}
|z|\le \frac{C^+_z}{\sqrt{1+t}}\ , \label{eq4.5}
\end{equation}
where $C_z^+:=C^+_{|\nabla u |}$.

\bigskip
Next, from (\ref{eq1.10}) it is easy to derive (see Ch 2 of
\cite{List1}) that
\begin{equation}
\frac{\partial R}{\partial t} = \Delta R +2R_{ij}R^{ij}
+2k^2_n (\Delta u)^2 -2k^2_n |\nabla \nabla u|^2
-4k^2_n R_{ij}\nabla^i u \nabla^j u+\nabla_X R\ . \label{eq4.6}
\end{equation}
Unlike in the Ricci flow, this equation does not lead to the
preservation of scalar curvature. However, defining
\begin{equation}
S_{ij}:=R_{ij}-k^2_n \nabla_i u \nabla_j u\ , \label{eq4.7}
\end{equation}
and using (\ref{eq1.12}), then (\ref{eq4.1}) leads to
\begin{eqnarray}
\frac{\partial S}{\partial t}&=&\Delta S +2S_{ij}S^{ij}
+2k^2_n (\Delta u)^2 +\nabla_X S \nonumber\\
&=& \Delta S + 2 \left ( S_{ij} -\frac1n g_{ij}S \right )
\left ( S^{ij}-\frac1n g^{ij}S \right ) +\frac2n S^2
+2k^n_2 (\Delta u )^2+\nabla_X S \nonumber\\
&\ge& \Delta S + \frac2n S^2+\nabla_X S\ , \label{eq4.8}
\end{eqnarray}
where $S:=g^{ij}S_{ij}$. Then we obtain

%PROPOSITION 4.4%%%%%%%%%%%%%%%%%%%%%%%%%%%%%%%%%%%%%%%%%
\bigskip
\noindent{\bf Proposition 4.4.} {\sl Let $(g,u)$ be an
asymptotically flat solution of the flow on the domain $[0,T]\times
{\mathbb R}^n$, $0<T<{\tilde T}$. For any $n\ge 2$
\begin{equation}
S :=R-k^2_n|\nabla u|^2\ge \frac{C^-_S}{1+t}\ , \label{eq4.9}
\end{equation}
where $C^-_S \le 0$ is a constant depending only on the initial data
and $n$, and if $S(0,x)\ge 0$ for all $x\in{\mathbb R}^n$, then
$S(t,x)\ge 0$ for all $t\ge 0$ and all $x\in{\mathbb R}^n$.}

\bigskip

\noindent{\bf Proof.} Use (\ref{eq4.8}) to apply Lemma 4.1 to $-S$.
The $k^2$ of Lemma 4.1 takes the value $-2/n$. This yields
$C^-_S=\min \left \{ -\frac{n}{2}, \inf_x S(0,x) \right \}$.\qed

\bigskip

Note that List's flow does not necessarily preserve positive scalar
curvature, though it does preserve the positivity of $R-k^2_n
|\nabla u |^2$.

%Rot. symm. bounds%%%%%%%%%%%%%%%%%%%%%%%%%%%%%%%%%%%%%%%%
\subsection{Bounds that hold in rotational symmetry}

\noindent The bounds of the previous subsection are valid with or
without rotational symmetry. We now specialize to the rotationally
symmetric flow equations (\ref{eq3.11}, \ref{eq3.12}) on
$[0,T]\times [0,\infty)\ni (t,r)$, $0<T<{\tilde T}$.

%Bounds on f%%%%%%%%%%%%%%%%%%%%%%%%%%%%%%%%%%%%%%%%%%%%
\subsubsection{Bounds on $f$}

\noindent In this subsection, we derive bounds on
$f:=\sqrt{g_{rr}}$. These bounds allow us to address two concerns.
The first is that our coordinate system may break down during the
flow. This will happen if $f$ diverges to $+\infty$ or approaches
zero along the flow. Note that $f$ diverges at some $r>0$ iff the
mean curvature $H$ of the $r=const$ sphere goes to zero. The mean
curvature is given by
\begin{equation}
H=\frac{n-1}{rf}\ . \label{eq4.10}
\end{equation}
Thus, divergence of $f$ at finite $r$ implies the presence or
formation of a minimal hypersphere. We will show that this cannot
happen.

The second concern arises because positivity of the scalar curvature
is not strictly preserved along the flow, even though the results of
the previous subsection show that $R$ is bounded below and the bound
tends to zero in time. One may wonder whether this is enough for
purposes of the static minimization conjecture, where one seems to
want the static metric to arise as a limit of a sequence of positive
scalar curvature metrics.

Here we show that the rotationally symmetric flow (\ref{eq3.11},
\ref{eq3.12}) on ${\mathbb R}^3$ does preserve the positivity of the
Brown-York and Misner-Sharpe quasi-local masses. The {\it Brown-York
mass} of a closed embedded hypersurface
$\Sigma\hookrightarrow{\mathbb R}^3$ is defined to be $\mu_{\rm
BY}[\Sigma]:=\int_{\Sigma} \left ( H_0 -H \right ) d\Sigma$, where
$H$ is the mean curvature of $\Sigma$ and $H_0$ is the mean
curvature of an isometrically embedded image of $\Sigma$ in flat
space. In our case, for a sphere of radius $r$ about the origin, we
have
\begin{equation}
\mu_{\rm BY} [\Sigma] := \frac{8\pi}{r} \left (
1-\frac{1}{f(t,r)} \right ) \ . \label{eq4.11}
\end{equation}
The {\it Misner-Sharpe mass} is defined only for rotationally
symmetric metrics and is given for $n=3$ by
\begin{equation}
\mu_{\rm MS}:=\frac{8\pi}{r} \left ( 1-\frac{1}{f^2(t,r)} \right )
= \left ( 1+\frac{1}{f}
\right )\mu_{\rm BY}\ , \label{eq4.12}
\end{equation}
so it is positive if and only if the Brown-York mass is.

We now show that for any finite $t$ along the flow, $f(t,r)$ is
bounded above and below. As a result, if no minimal sphere is
present initially then none will form, and for $n=3$ if $\mu_{\rm
BY}$ is initially positive then it will always be so (likewise for
$\mu_{\rm MS}$). In fact, this will hold in any dimension if we take
(\ref{eq4.11}) and (\ref{eq4.12}), without modification, to be the
definitions of $\mu_{\rm BY}$ and $\mu_{\rm MS}$ in any dimension
(this is not what is usually done, however).

\bigskip

%PROPOSITION 4.5%%%%%%%%%%%%%%%%%%%%%%%%%%%%%%%%%%%%%%%%%
\noindent {\bf Proposition 4.5.}  {\sl
\begin{equation}
C^-_f\le f(t,r)\le C^+_f(1+t)^p\ ,  \label{eq4.13}
\end{equation}
where $p=1+\left ( k_nC^+_z\right )^2$ for all $n\ge 2$. The
constants $C^{\pm}_f$ depend only on the initial data $\{
(f(0,r),z(0,r))\}$ and, for $n=2$, $f_{\infty}$ (equivalently, the
$n=2$ mass).}

\bigskip
\noindent {\bf Proof.} Let $w(t,r)=f^2(t,r)-1$. Then (\ref{eq3.11})
yields
\begin{equation}
\frac{\partial w}{\partial t}
=\frac{1}{f^2}\frac{\partial^2 w}{\partial r^2}
-\frac{3}{2f^4} \left (\frac{\partial w}{\partial r} \right )^2
+\left (\frac{n-2}{r}-\frac{1}{rf^2} \right )
\frac{\partial w}{\partial r}-\frac{2(n-2)}{r^2}w
+2k^2_nf^2z^2\ , \label{eq4.14}
\end{equation}
subject to the boundary conditions
\begin{equation}
w(t,r)\to
\begin{cases}
0 &\text{for $r\to 0$ and all $n\ge 2$,}\\
w_{\infty}:=f^2_\infty-1 &\text{for $r\to\infty$ and $n=2$, and} \\
0\ ,&\text{for $r\to\infty$ and $n>2$.}
\end{cases}\label{eq4.15}
\end{equation}

Now consider the closed annular region $A_{\epsilon}(T):=
[0,T]\times \left [ \epsilon,\frac{1}{\epsilon}\right ] \ni (t,r)$.

\begin{enumerate}
\item[{\it (i)}] {\it Case $n>2$:} By (\ref{eq4.15}), $\inf \{
    w(t,r)| 0\le t \le T, r\ge 0\} \le 0$. We observe from
    (\ref{eq4.14}) that if $w(t,r)$ has a negative minimum in
    $A_\epsilon(T)$, such a minimum must lie on the parabolic
    boundary.\footnote
{e.g., Let $\Psi=-w$ and then observe that the inequality
(\ref{eq4.1}) applies on $A_{\epsilon}(T)$ (with $k=0$), so we
can use Lemma 4.1(i).}
    Taking $\epsilon$ sufficiently large, then by the
    boundary conditions, the negative minimum must lie on the
    initial boundary, so then $w(t,r)\ge \inf_r \left\{ w(0,r)
    \right\}$.
\item[{\it (ii)}] {\it Case $n=2$:} We consider the function
    $W_{\epsilon}(t,r)=w(t,r)+\epsilon\cdot t$ on
    $A_{\epsilon}(T)$. From (\ref{eq4.14}), $W$ has the
    following evolution equation:
\begin{equation}
\frac{\partial W_{\epsilon}}{\partial t}
=\frac{1}{f^2}\frac{\partial^2 W_{\epsilon}}{\partial r^2}
-\frac{3}{2f^4}\left(\frac{\partial W_{\epsilon}}{\partial t}\right)^2
-\frac{1}{rf^2}\frac{\partial W_{\epsilon}}{\partial r} +2k^2_nf^2z^2
+\epsilon\ , \label{eq4.16}
\end{equation}
with boundary conditions $W_{\epsilon}(t,r)\rightarrow
\epsilon\cdot t$ as $r \rightarrow 0$ and
$W_{\epsilon}(t,r)\rightarrow w_\infty+\epsilon\cdot t$ as $r
\rightarrow \infty$. For $\delta
>0$, the minimum of $W_{\epsilon}$ on $A_\epsilon(T)$ must lie on the
parabolic boundary. Taking $\epsilon\to 0$, then we obtain
$w(t,r) \ge \min \left\{0,w_\infty, \inf_r\left \{ w(0,r)\right
\} \right\}$.
\end{enumerate}

\noindent This proves the left-hand (i.e., inferior) inequality in
(\ref{eq4.13}) and shows that
\begin{equation}
C^-_f=
\begin{cases}
\inf_r \left\{ f(0,r) \right\}&\text{for $n>2$,}\\
\min\left \{ 1, f_{\infty}, \inf_r\left \{ f(0,r)\right \}
\right\}&\text{for $n=2$.}
\end{cases}\label{eq4.17}
\end{equation}

To prove the superior inequality, consider the function
\begin{equation}
Q(t,r):=\frac{w(t,r)}{(1+t)^{2p}}\ . \label{eq4.18}
\end{equation}
For $t<\tau$, $Q$ obeys $Q(0,r)=w(0,r)$, $Q(t,0)=0$, and either
$\lim_{r\to\infty}Q(t,r)=0$ if $n>2$ or
$\lim_{r\to\infty}Q(t,r)=\frac{w_{\infty}}{(1+t)^{2p}}$ if $n=2$. As
well, $Q$ solves the PDE
\begin{eqnarray}
\frac{\partial Q}{\partial t}
&=&\frac{1}{f^2}\frac{\partial^2Q}{\partial r^2}
-\frac{3(1+t)^{2p}}{2f^4}\left ( \frac{\partial Q}{\partial r}
\right )^2 +\left ( \frac{n-2}{r}-\frac{1}{rf^2}\right )
\frac{\partial Q}{\partial r}\nonumber \\
&&+\left ( 2k_n^2z^2-\frac{2(n-2)}{r^2}-\frac{{2p}}{1+t}\right )
Q +\frac{2k_n^2}{(1+t)^{2p}}z^2 \nonumber\\
&\le&\frac{1}{f^2}\frac{\partial^2Q}{\partial r^2}
-\frac{3(1+t)^{2p}}{2f^4}\left ( \frac{\partial Q}{\partial r}
\right )^2 +\left ( \frac{n-2}{r}-\frac{1}{rf^2}\right )
\frac{\partial Q}{\partial r}\nonumber \\
&&+\left ( \frac{2k_n^2(C^+_z)^2}{1+t}
-\frac{2(n-2)}{r^2}-\frac{{2p}}{1+t}\right ) Q
+\frac{2k_n^2(C^+_z)^2}{(1+t)^{{2p}+1}} \label{eq4.19}
\end{eqnarray}
using (\ref{eq4.5}). Choose
\begin{equation}
p:=1+\left ( k_nC^+_z\right )^2\ . \label{eq4.20}
\end{equation}
Then (\ref{eq4.19}) yields
\begin{equation}
\frac{\partial Q}{\partial t}
\le\frac{1}{f^2}\frac{\partial^2Q}{\partial r^2}
+\left ( \frac{n-2}{r}-\frac{1}{rf^2}\right )
\frac{\partial Q}{\partial r} +\frac{2}{1+t}
\left ( \frac{(k_nC^+_z)^2}{(1+t)^{2p}}-Q\right )\label{eq4.21}
\end{equation}
whenever $Q\ge 0$, where we've discarded some manifestly negative
terms. Clearly this equation does not permit $Q$ to have a maximum
on the parabolic interior of $[0,T]\times [\epsilon,1/\epsilon]$
unless $Q\le \frac{ (k_n C^+_z)^2}{(1+t)^{2p}}$, whence by
(\ref{eq4.18}) we get $w(t,r)\le \left ( k_nC^+_z\right )^2$ and
then
\begin{equation}
f(t,r)\le \sqrt{1+(k_n C^+_z)^2}\ . \label{eq4.22}
\end{equation}
Otherwise, the maximum of $Q$ can occur on the parabolic boundary.
Then taking $\epsilon$ sufficiently small, if a positive maximum for
$Q$ occurs either:

\begin{enumerate}
\item[{\it (a)}] The maximum of $Q$ occurs on the initial
    boundary $t=0$. This can occur for any $n\ge 2$. Using
    (\ref{eq4.18}), $w=f^2-1$, and the fact that
    $Q(0,r)=w(0,r)=f^2(0,r)-1$, then $f^2(t,r)\le
    1+(1+t)^{2p}\sup_r\{ f^2(0,r)-1 \}\le (1+t)^{2p}\sup_r\{
    f^2(0,r) \}$. Combining this with (\ref{eq4.22}) yields
\begin{equation}
f(t,r) \le \max \left \{\sqrt{1+(k_n C^+_z)^2},
(1+t)^p\sup_r\left \{ f(0,r) \right \} \right \}\le C^+_f (1+t)^p
\label{eq4.23}
\end{equation}
for $C^+_f=\max \left \{\sqrt{1+(k_n C^+_z)^2}, \sup_r\left \{
f(0,r) \right \} \right \}$, or
\item[{\it (b)}] $n=2$ and the maximum of $Q$ is
    $\frac{w_{\infty}}{(1+t)^{2p}}< \max\{0,w_{\infty} \}$.
    Combining this with (\ref{eq4.22}) and (\ref{eq4.23}), we
    obtain
\begin{equation}
f(t,r)\le \max \left \{ f_{\infty},\sqrt{1+(k_n C^+_z)^2},(1+t)^p\sup_r
\left \{ f(0,r) \right \} \right \}\le C^+_f (1+t)^p\ , \label{eq4.24}
\end{equation}
\end{enumerate}
for $C^+_f=\max \left \{f_{\infty}, \sqrt{1+(k_n C^+_z)^2},
\sup_r\left \{ f(0,r) \right \} \right \}$.\qed

%COROLLARY 4.6%%%%%%%%%%%%%%%%%%%%%%%%%%%%%%%%%%%%%%%%%
\bigskip
\noindent{\bf Corollary 4.6.} {\sl (i) If no minimal hypersphere is
present initially, none forms at any $t<\infty$. (ii) For $n=2$, if
the Brown-York mass $\mu_{\rm BY}(0,r)$ of every $r=const$
hypersphere about the origin is $\ge 0$ at $t=0$, then $\mu_{\rm
BY}(t,r)\ge 0$ for every $r\in{\mathbb R}$ and every $t>0$; the same
holds for the Misner-Sharpe mass.}

\bigskip
\noindent{\bf Proof.} The first statement follows immediately from
(\ref{eq4.23}) and (\ref{eq4.10}). The second statement follows from
(\ref{eq4.17}), (\ref{eq4.13}), and (\ref{eq4.11}) (or, for the
Misner-Sharpe mass, ({\ref{eq4.12})). \qed

\bigskip
\noindent{\bf Remark 4.7.} If the assumptions of (ii) hold and if
the flow converges $({\mathbb R}^3, g(t), 0)$ in the pointed
Cheeger-Gromov sense to $({\mathbb R}^3,g_{\infty},0)$, then
$({\mathbb R}^3, g_{\infty})$ will have nonnegative Brown-York mass
at each $r$ (by (\ref{eq4.11}) and the fact that the sign of
$1-\frac{1}{f}$ will be preserved under the diffeomorphisms
$({\mathbb R}^3, g(t), 0)\to ({\mathbb R}^3,g_{\infty},0)$). Since
$\lim_r \mu_{\rm BY}=m_{\rm ADM}:=\frac{1}{16\pi}
\int_{S^2_{\infty}}\delta^{ij} \left ( g_{ki,j}-g_{ij,k}\right )
dS^i$ (we take this limit along $r=const\to\infty$ spheres), the ADM
mass of the limit manifold will be nonnegative.

\subsubsection{Proof of Proposition 3.1}

\noindent{\bf Proof.} Setting $t=0$ in (\ref{eq4.13}), we see that
$C^{-}_f \le f(0,r)\equiv a(r)\equiv q\left (0,\rho(r)\right ) \le
C^+_f$. Assume, by way of contradiction, that ${\tilde T}<T_M$.
Then, by Proposition 2.4, there are constants $K$ and $C$ such that
\begin{equation}
e^{-(2nK+4C)T}\left ( C^-_f \right )^2 \leq q^2(t,\rho)
\leq e^{(2nK+4C)T}\left ( C^+_f \right )^2 \label{eq4.25}
\end{equation}
for $0\le t \le T \le {\tilde T}$. Furthermore, (\ref{eq4.13}) holds for
all $t\in [0,T]$ so, dividing (\ref{eq4.25}) by
(\ref{eq4.13}) and using (\ref{eq3.6}), we get
\begin{equation}
e^{-(2nK+4C)T}\left ( \frac{C^-_f}{C^+_f(1+T)^p} \right )^2 \leq
\frac{q^2(t,\rho(r))}{f^2(t,r)}=\left ( \frac{\partial h}{\partial \rho} \right )^2
\leq e^{(2nK+4C)T}\left ( \frac{C^+_f}{C^-_f} \right )^2 \label{eq4.26}
\end{equation}
for all $t\in [0,T]$. We can replace $T$ using
that $T\le {\tilde T}<T_M$, obtaining
\begin{equation}
e^{-(2nK+4C)T_M}\left ( \frac{C^-_f}{C^+_f(1+T_M)^p} \right )^2 \leq
\left ( \frac{\partial h}{\partial \rho} \right )^2
\leq e^{(2nK+4C)T_M}\left ( \frac{C^+_f}{C^-_f} \right )^2 \ . \label{eq4.27}
\end{equation}
By comparison, we see that the constants in the inequality
(\ref{eq3.4}) are in fact independent of $T$. Since the inequalities
hold for any $T<{\tilde T}$ and are closed relations, they hold for
$T={\tilde T}$ as well and, by adjusting the constants slightly if
necessary (keeping the inferior one positive of course), then
(\ref{eq3.4}) holds for $t$-values beyond ${\tilde T}$,
contradicting the assumption. \qed

%Sectional curvature bound%%%%%%%%%%%%%%%%%%%%%%%%%%%%%%%%%%%%%
\bigskip
\subsubsection{A bound on tangential sectional curvature}

\noindent We will now obtain a bound on the behaviour of $f$ at the
origin. This is in fact a lower bound on $\lambda_2$, which for
$n\ge 3$ is the sectional curvature in planes tangent to the
$r=const$ spheres.

%PROPOSITION 4.8%%%%%%%%%%%%%%%%%%%%%%%%%%%%%%%%%%%%%%%%%
\bigskip
\noindent {\bf Proposition 4.8.}  {\sl For all $n\ge 2$,
$\lambda_2(t,r)$ is bounded below by a constant $C^-_{\lambda_2}>0$
which depends only on the initial data $f(0,r)$ such that
$\lambda_2(t,r)\ge-C^-_{\lambda_2}/(1+t)$.}

\bigskip
\noindent {\bf Proof.} In close (but not exact) analogy to
\cite{OW}, we will approximate $\lambda_2$ by a sequence of
functions $u_m(t,r)$, $0<m<2$, defined by
\begin{eqnarray}
u_m(t,r)&:=&\left(\frac{2}{r^m+r^2}\right)\left ( 1-\frac{1}{f^2}
\right ) \mathrm{\ for\ } r>0\ ,
\label{eq4.28}\\
u_m(t,0)&:=&\lim_{r\rightarrow0} u_m(t,r)\ .
\label{eq4.29}
\end{eqnarray}
The $u_m(t,r)$ functions have the following useful properties:
\begin{enumerate}
\item[\it (i)] $u_m(t,0)=0$ for all $0<m<2$ and
    $\lim_{r\rightarrow \infty}u_m(t,r)=0$ for all $0<m\le2$.
\item[{\it (ii)}] For fixed $t$ and $r\ne0$, the map $m \mapsto
    u_m(t,r)$ is continuous at $m=2$, and in fact
\begin{equation}
\lambda_2=u_2\ .
\label{eq4.30}
\end{equation}
\end{enumerate}
Now define new functions
\begin{equation}
U_m(t,r)=(1+t)u_m(t,r)\ , \label{eq4.31}
\end{equation}
and note that $U_m(0,r)=u_m(0,r)$. From (\ref{eq3.11}) we obtain an
evolution equation for $U_m(t,r)$ given by
\begin{eqnarray}
\frac{\partial U_m}{\partial t}
&=&\frac{1}{f^2}\frac{\partial^2U_m}{\partial r^2}
+\frac{(r^m+r^2)}{4(1+t)}\left(\frac{\partial U_m}{\partial r}\right)^2
+\frac{(2r+mr^{m-1})}{2(1+t)}U_m\frac{\partial U_m}{\partial r} \nonumber\\
&&+\left [\frac{2(2r+mr^{m-1})}{f^2(r^m+r^2)} -\frac{1}{rf^2}
+\frac{(n-2)}{r}\right]\frac{\partial U_m}{\partial r}\nonumber \\
&&+\frac{1}{2(1+r^{2-m})(1+t)} \bigg [ (4-m)(m+n-2)+m(n-2)+2(n-1)r^{2-m}
\nonumber\\
&&+r^{m-2}(2-m)(m+n-2)+r^{m-2}m\left(\frac{m}{2} +n-2\right)\bigg ]
U_m^2 \nonumber \\
&&-\frac{(2-m)(m+n-2)}{r^2(1+r^{2-m})}U_m
+\left(\frac{4(1+t)}{r^m+r^2}\right )k_n^2z^2 \nonumber \\
&&-2k_n^2z^2U_m +\frac{1}{(1+t)}U_m\label{eq4.32}\\
&\ge&\frac{1}{f^2}\frac{\partial^2U_m}{\partial r^2}
+\left [ \frac{(2r+mr^{m-1})}{2(1+t)}U_m
+\frac{2(2r+mr^{m-1})}{f^2(r^m+r^2)} -\frac{1}{rf^2}
+\frac{(n-2)}{r}\right]
\frac{\partial U_m}{\partial r} \nonumber\\
&&+\frac{1}{1+t}\left [ (n-1)U_m^2+U_m \right ]
-\frac{(2-m)(m+n-2)}{r^2(1+r^{2-m})}U_m-2k_n^2z^2U_m
\ , \label{eq4.33}
\end{eqnarray}
where the inequality holds at least for $1\le m <2$ and $n\ge 2$.
Furthermore, if
\begin{equation}
U_m<-\frac{1}{n-1}\ , \label{eq4.34}
\end{equation}
we then obtain
\begin{eqnarray}
\frac{\partial U_m}{\partial t}&>&\frac{1}{f^2}
\frac{\partial^2U_m}{\partial r^2}
+\bigg [ \frac{(2r+mr^{m-1})}{2(1+t)}U_m
+\frac{2(2r+mr^{m-1})}{f^2(r^m+r^2)}\nonumber\\
&&\qquad -\frac{1}{rf^2}
+\frac{(n-2)}{r}\bigg ]
\frac{\partial U_m}{\partial r}\label{eq4.35}
\end{eqnarray}

As with Proposition 4.5, we work first on the annulus
$A_{\epsilon}(T)$. From (\ref{eq4.34}) and (\ref{eq4.35}), we see
that $U_m$ cannot have a minimum $<-\frac{1}{n-1}$ at some
$(t_0,r_0)$ in the parabolic interior of the annulus. The minimum,
if $<-\frac{1}{n-1}$, must lie on the parabolic boundary of
$A_{\epsilon}(T)$. Taking $\epsilon\to 0$ and recalling that
$u_m(t,\epsilon)\to 0$ and $u_m(t,1/\epsilon)\to 0$, whence
$U_m(t,\epsilon)\to 0$ and $U_m(t,1/\epsilon)\to 0$ as well, then
the minimum, if $<-\frac{1}{n-1}$, must lie at $t=0$; that is,
\begin{eqnarray}
U_m(t,r)&\ge& \min \left \{ -\frac{1}{n-1}, \inf_r \left\{U_m(0,r)
\right \} \right \}=\min \left \{ -\frac{1}{n-1}, \inf_r \left
\{ u_m(0,r) \right \} \right \}
\nonumber\\
&=& \min \left \{ -\frac{1}{n-1}, \inf_r\left \{ \frac{2}{r^m+r^2}
\left ( 1-\frac{1}{f^2(0,t)}\right )\right \}\right \} \nonumber\\
&\ge& \min \left \{ -\frac{1}{n-1},
\inf_r\left \{ \frac{2}{r^2}\left ( 1-\frac{1}{f^2(0,t)}\right )
\right \} \right \}\nonumber\\
&=& \min \left \{ -\frac{1}{n-1}, 2\inf_r \{ \lambda_2(0,r) \} \right \}
=:-C_{\lambda_2}^-\ , \label{eq4.36}
\end{eqnarray}
where $C_{\lambda_2}^-\ge 0$. We now take $m\nearrow 2$ in
(\ref{eq4.36}), so that $U_m\to (1+t)\lambda_2$ by (\ref{eq4.30})
and (\ref{eq4.31}). Using (\ref{eq4.34}) and (\ref{eq4.35}) as well,
(\ref{eq4.36}) yields
\begin{equation}
\lambda_2(t,r)\ge -\frac{C^-_{\lambda_2}}{(1+t)} \ , \label{eq4.37}
\end{equation}
\qed

%Smoothness of |\nabla u|%%%%%%%%%%%%%%%%%%%%%%%%%%%%%%%%%%%%%%
\bigskip
\subsubsection{Smoothness of $|\nabla u|$ for $n=2$}

%PROPOSITION 4.9%%%%%%%%%%%%%%%%%%%%%%%%%%%%%%%%%%%%%%%%%
\bigskip
\noindent {\bf Proposition 4.9.} {\sl Assume $n=2$. Then
\begin{equation}
\frac{1}{r} \left |\nabla u (t,r)\right \vert \le C^+_{\zeta}
\label{eq4.38}\ ,
\end{equation}
where the constant $C^+_{\zeta}$ depends only on the (smooth)
initial data for $\nabla u$.}

\bigskip
\noindent {\bf Proof.} Let $\zeta_m(t,r)=2\frac{z(t,r)}{r+r^m}$ for
$0<m<1$. Computing from (\ref{eq3.12}) we then obtain that $\zeta_m$
obeys
\begin{eqnarray}
\frac{\partial \zeta_m}{\partial t}&=&\frac{1}{f^2}
\frac{\partial^2 \zeta_m}{\partial r^2}
+\left\{\frac{3+(2m+1)r^{m-1}}{rf^2(1+r^{m-1})}
+\frac{(n-2)}{r} \right\}\frac{\partial \zeta_m}{\partial r}\nonumber \\
&&+\frac{(m-1)r^{m-1}}{r^2(1+r^{m-1})}
\left\{\frac{m+1}{f^2} +(n-2)\right\}\zeta_m \nonumber \\
&&+(n-2)\lambda_2\zeta_m-k^2_nz^2\zeta_m.
\label{eq4.39}
\end{eqnarray}

For $\zeta_m>0$, $m<1$, and $n=2$, (\ref{eq4.39}) reduces to
\begin{equation}
\frac{\partial \zeta_m}{\partial t}\le \frac{1}{f^2}
\frac{\partial^2 \zeta_m}{\partial r^2}
+\left [ \frac{3+(2m+1)r^{m-1}}{rf^2(1+r^{m-1})}
\right ] \frac{\partial \zeta_m}{\partial r}\ .
\label{eq4.40}
\end{equation}
As usual, restrict attention to the annulus
$A_{\epsilon}(T):=[0,T]\times [\epsilon, 1/\epsilon]$, for some
chosen $\epsilon>0$ and $T<\tau$, $\tau$ as above. By the maximum
principle, $\zeta_m$ must have a maximum in $A_{\epsilon}(T)$, but
by (\ref{eq4.40}) this cannot occur in the parabolic interior of
$A_{\epsilon}(T)$. If the maximum occurs at $r=\epsilon$, then take
$\epsilon\to 0$. By regularity, $\frac{\partial u}{\partial r}\in
\mathcal{O}(r)$ as $r \rightarrow 0$, so
$\zeta_m(t,\epsilon)\in{\cal O}(\epsilon^{1-m})\to 0$, $\forall t\in
[0,T]$. Similarly, at large $r$, $\frac{\partial u}{\partial r}\in
\mathcal{O}(1/r)$ and so $\zeta_m\to 0$ as $r=1/\epsilon \to
\infty$. Thus, for $\epsilon$ small enough, $\zeta_m$ cannot have a
positive maximum at any $t>0$, and since the supremum of $\zeta_m$
is nonnegative on the $t=0$ boundary then
\begin{equation}
\left |\zeta_m(t,r)\right \vert \le \sup_r\{\left \vert
\zeta_m(0,r)\right \vert \}\le \sup_r\left \{
\frac{2}{r}|z(0,r)|\right \} =:C^+_{\zeta}\ , \label{eq4.41}
\end{equation}
where the smoothness of the initial data is used to infer the
boundedness of $|z(0,r)|/r$. Finally, since $C^+_{\zeta}$ is
independent of $m$, we can take $m \nearrow 1$ to complete the
argument.\qed

\bigskip

%Summary of a priori bounds%%%%%%%%%%%%%%%%%%%%%%%%%%%%%%%%%%%%
\subsection{Summary of {\it a priori} bounds}

\noindent In summary, we have the following bounds for all $t\in
[0,T]$, for all $x\in{\mathbb R}^n$, and, assuming rotational
symmetry, for all $r\in[0,\infty)$.

\begin{enumerate}
\item $const\le f^2 \le const\cdot (1+t)^p$.
\item $R(t,x) \ge -\frac{const}{1+t}$.
\item $R(t,x)\ge k^2_n |\nabla u|^2$ for all $(t,x)$ if it holds
    at $t=0$.
\item $\lambda_2(t,r)\ge -\frac{const}{(1+t)}$ and
    $\lambda_2(t,r)\ge 0$ if $f(0,r)\ge 1$ for all $r$.
\item $|\nabla u|^2\le \frac{const}{1+t}$.
\item $\frac{1}{r}|\nabla u|\le const$ when $n=2$.
\end{enumerate}
The constants denoted $const$ here are positive and distinct.
These constants, and $p$, depend only on the initial data, $n$,
$k_n$, and (for $n=2$) $f_{\infty}$, and do not depend on $T$, as
can be seen by the explicit expressions for the constants given in
the preceding section.

For $n\ge 3$, we can summarize the picture that these bounds present
as follows.

\bigskip
\noindent{\bf Proposition 4.10.} {\sl Assume $n\ge 3$. One of the
following possibilities holds:
\begin{enumerate}
\item The flow (\ref{eq3.12}, \ref{eq3.13}) exists for all
    $(t,r)\in [0,\infty) \times [0,\infty)$.
\item There is a sequence of points $(t_k,r_k)$ with $r_k\to 0$
    such that $\lambda_1(t_k,r_k) = \lambda_2(t_k,r_k)\to
    +\infty$ as $t_k\nearrow T_M$.
\item There is a sequence of points $(t_k,r_k)$ along which
    $\lambda_1(t_k,r_k) \to +\infty$ as $t_k\nearrow T_M$ but
    $\lambda_2(t_k,r_k)$ remains bounded along every such
    sequence.
\end{enumerate}
}
\bigskip

\noindent{\bf Remark 4.11.} The considerations of the next section
will eliminate the third possibility from this list.
\bigskip

\noindent{\bf Corollary 4.12.} {\sl Either the flow exists for all
$t>0$ or $\limsup_{t\nearrow T_M} \sup_r \lambda_1 = \infty$ and
$\limsup_{t\nearrow T_M} \sup_r R = \infty$.}
\bigskip

\noindent{\bf Remark 4.13.} The Corollary is also true for $n=2$,
since then $2\lambda_1=R\ge -\frac{const}{1+t}$.
\bigskip

\noindent{\bf Proof of 4.10.} By the continuation principle, either
the flow exists for all $t>0$ or at least one sectional curvature
diverges as $t\to T_M$. We first consider the case of
$\lambda_2\to\infty$. Then there is a sequence of points
$(t_k,r_k)$, $t_k<t_{k+1}<T_M$, along which $\lambda_2$ assumes
successive maximum values; $\lambda_2(t_k,r_k)\ge \lambda_2(t,r)$
for all $(t,r)\in [0,t_k]\times [0,\infty)$. From the definition of
$\lambda_2$, we see that $r_k\to 0$ along any such sequence, and at
these points, the Bianchi identity (\ref{eq3.16}) shows that
$\lambda_1(t_k,r_k)=\lambda_2(t_k,r_k)$.

Since $\lambda_2\ge -\frac{const}{1+t}$, the only remaining cases
are those for which $\lambda_2$ remains bounded. But then
$\lambda_1$ cannot diverge to $-\infty$ because $(n-1)\left (
2\lambda_1+(n-2)\lambda_2\right ) \equiv R \ge -\frac{const}{1+t}$.
Thus we have eliminated all possibilities that are not enumerated in
the proposition. \qed
\bigskip

%SECTION 5%%%%%%%%%%%%%%%%%%%%%%%%%%%%%%%%%%%%%%%%%%%%
\section{Smoothness of $|\nabla u|$ and long-time existence}
\setcounter{equation}{0}

\subsection{Smoothness of $|\nabla u|$ and an upper bound on $\lambda_2$}

\noindent In this section, we will show that whenever
$\frac{1}{r}|\nabla u|$ remains finite, the flow exists for all
$t>0$. When $n=2$, we have already shown in Proposition 4.9 that
$\frac{1}{r}|\nabla u|$ remains bounded. In this subsection, the
first proposition we present shows that this will also be the case
for $n>3$, provided that $\lambda_2$ remains finite. In consequence,
the $n\ge 3$ flow will fail to exist only if $\lambda_2$ diverges at
finite $T$. We then show that, conversely, when $\frac{1}{r}|\nabla
u|$ remains finite, so does $\lambda_2$. This follows for all
dimensions, including $n=2$, and is useful when $n=2$ even though
the combination $\lambda_2=\frac{1}{r^2}\left (1 -\frac{1}{f^2}
\right)$ is of course not a sectional curvature in that case.

\bigskip
\noindent{\bf Proposition 5.1.} {\sl Assume that there is a function
$F_{\lambda_2}^+:[0,\infty)\to [0,\infty)$ such that
$\lambda_2(t,r)\le F_{\lambda_2}^+(T)$ for all $0\le t\le T$. Then
there is a function $F_{\zeta}^+:[0,\infty)\to [0,\infty)$ such that
\begin{equation}
\frac{z}{r}\equiv\frac{1}{r} \left |\nabla u (t,r)\right \vert
\le F^+_{\zeta}(T) \ {\rm whenever\ } 0\le t \le T \ .
\label{eq5.1}
\end{equation}
}

\bigskip

\noindent{\bf Proof.} By Proposition 4.9, this is true for $n=2$
(without the assumption on $\lambda_2$ and with $F^+_{\zeta}=
C^+_{\zeta}=const$). Thus, assume $n\ge 3$. Choose some $T>0$ and
define $\xi_m:= \zeta_m/(1+t)^{(n-2)F^+_{\lambda_2}(T)}$ for $0\le t
\le T$. Then from (\ref{eq4.39}) with $0<m<1$ and $\xi_m>0$, we
obtain
\begin{equation}
\frac{\partial \xi_m}{\partial t}\le\frac{1}{f^2}
\frac{\partial^2 \xi_m}{\partial r^2}
+\left [ \frac{3+(2m+1)r^{m-1}}{rf^2(1+r^{m-1})}
\right ] \frac{\partial \xi_m}{\partial r}\ , \label{eq5.2}
\end{equation}
and the proof proceeds precisely as in Proposition 4.9. This implies
that $\xi_m$ is bounded above by a constant depending only on
initial data, and thus $\zeta_m\le F^+_{\zeta}(T) :=const \cdot
(1+t)^{(n-2)F^+_{\lambda_2}(T)}$. Since $F^+_{\lambda_2}$ is defined
for all $T>0$, so is $F^+_{\zeta}$, and since the bound is
$m$-independent, we extend to $m=1$. \qed

%PROPOSITION 5.2%%%%%%%%%%%%%%%%%%%%%%%%%%%%%%%%%%%%%%%%%

\bigskip
\noindent {\bf Proposition 5.2.}  {\sl Conversely, assume that
equation (\ref{eq5.1}) holds for all $T>0$. Then there is a function
$F^+_{\lambda_2}:[0,\infty)\times [0,\infty)$ such that}
\begin{equation}
\lambda_2(t,r)\le F^+_{\lambda_2}(T)\ {\rm whenever\ } 0\le t \le T
\ . \label{eq5.3}
\end{equation}

\bigskip
\noindent {\bf Proof.} We work, as always, on a compact annular
domain $A_{\epsilon}(T):=[0,T]\times [\epsilon,1/\epsilon]\ni
(t,r)$. Choose a positive function $F:[0,\infty)\to (0,\infty)$ such
that $F(0)=1$ and define functions
\begin{eqnarray}
V_m(t,r)&:=&\left(\frac{F(t)}{r^m+r^2}\right)\left( {f^2(t,r)}-1 \right)
\hspace{0.25cm}\mathrm{for}\hspace{0.25cm} r>0\ ,
\label{eq5.4}\\
V_m(t,0)&:=&\lim_{r\rightarrow0} V_m(t,r)\ .
\label{eq5.5}
\end{eqnarray}
Note that
\begin{equation} \lambda_2=\frac{2}{F(t)f^2(t,r)}V_2\ .
\label{eq5.6}
\end{equation}

We will show by the maximum principle that the $V_m(t,r)$ functions
have a uniform bound in $m$. From (\ref{eq3.11}), we obtain an
evolution equation for $V_m(t,r)$ given by
\begin{eqnarray}
\frac{\partial V_m}{\partial t}&=&\frac{1}{f^2}
\frac{\partial^2V_m}{\partial r^2}+\biggl [ \frac{2(mr^{m-1}+2r)}{f^2(r^m+r^2)}
-\frac{3(r^m+r^2)}{2f^4F}\frac{\partial V_m}{\partial r}\nonumber \\
&&
\qquad -\frac{3(mr^{m-1}+2r)}{f^4F}V_m +\frac{n-2}{r} -\frac{1}{rf^2}\biggr ]
\frac{\partial V_m}{\partial r}\nonumber \\
&&+\frac{r^{m-2}}{r^m+r^2}\left [ \frac{m(m-2)}{f^2}V_m+(n-2)(m-2)V_m+2k_n^2
\left \vert \nabla u \right \vert^2r^{2-m}F\right ] \nonumber\\
&&+\left [ 2k_n^2\left \vert \nabla u \right \vert^2+\frac{F^\prime}{F}
\right ] V_m-\frac{3(mr^{m-1}+2r)^2}{2f^4F(r^m+r^2)}V_m^2\ .
\label{eq5.7}
\end{eqnarray}

If $V_m(t,r)$ attains a positive maximum $V_m(t_0,r_0)<1$ for all
$m<2$, we are done, so assume to the contrary that the maximum is
$>1$. As well, for the moment assume that the maximum occurs at a
point $(t_0,r_0)$ in the parabolic interior of $A_{\epsilon}(T)$.
All the terms in (\ref{eq5.7}) that do not contain a derivative will
be negative provided that
\begin{equation}
\frac{2k_n^2\left \vert \nabla u \right \vert_0^2}{r_0^m+r_0^2}F(t_0)
+2k_n^2\left \vert \nabla u \right \vert_0^2V_m(t_0,r_0)
+\frac{F^\prime(t_0)}{F(t_0)}V_m(t_0,r_0)\le 0 \ .
\label{eq5.8}
\end{equation}
\noindent Observe that this implies that $F'(t_0)<0$, so take it to
be decreasing for all $t\ge 0$. Then $0<F(t)\le 1$ and so
(\ref{eq5.8}) will hold if it holds with $F(t_0)$ replaced by $1$ in
the first term on the left. Then, in the limit as $m\nearrow 2$,
this first term becomes $k^2_n\left(\frac{1}{r}\left \vert \nabla u
\right \vert_0\right)^2$. Since we assume that ({\ref{eq5.1}) holds,
we can control this term. Then we obtain the sufficient condition
\begin{eqnarray}
2k_n^2\left ( F^+_{\zeta}(T)\right )^2+\left [ 2k_n^2{(C^+_z)^2}
+\frac{F^\prime(t_0)}{F(t_0)}\right ] V_m(t_0,r_0)\le 0\ .
\label{eq5.9}
\end{eqnarray}

Since $V(t_0,r_0)\ge 1$ by assumption, a choice of $F$ that
satisfies this condition is
\begin{equation}
F(t)=e^{-Pt}\ , \ P:=2k_n^2\left [ \left (F^+_{\zeta}(T)\right )^2
+\left ( C^+_z \right )^2\right ] \ . \label{eq5.10}
\end{equation}

That is, choosing $F$ in (\ref{eq5.4}) to be given by
(\ref{eq5.10}), then either (i) $V_m$ is bounded above on
$A_{\epsilon}(T)$ by $1$ or the maximum of $V_m$ on
$A_{\epsilon}(T)$ resides on the parabolic boundary of
$A_{\epsilon}(T)$. On the spatial part of this boundary at
$r_0=1/\epsilon$, as $\epsilon\to 0$ we see from (\ref{eq5.4}),
(\ref{eq5.6}), and asymptotic flatness that $V_m\to 0$, so for any
fixed $\epsilon$ small enough, if the maximum were to occur on this
part of the boundary it would be less than $1$. Likewise, if it
occurs at $r_0=\epsilon$, then from (\ref{eq5.4}) we would have
$V_m(t_0,\epsilon)\sim F(t_0)(f^2-1)/\epsilon^m$, $m<2$, and then
local existence implies $V_m(t_0,\epsilon)\to 0$ as $\epsilon\to 0$,
so again for any fixed but sufficiently small $\epsilon$ we would
have $V_m<1$ at its maximum. Thus, the maximum, if $>1$, occurs on
the initial boundary, and so
\begin{eqnarray}
V_m(t,r)&\le&\max \left \{ 1, \sup_r \{ V_m(0,r)\} \right \}
\le \max \left \{ 1, \sup_r \{ 2V_2(0,r)\} \right \} \nonumber \\
&=& \max \left \{ 1, \sup_r \{ 2\lambda_2(0,r)\} \right \}  ,
\label{eq5.11}
\end{eqnarray}
for any $m<2$. Since the right-hand side is independent of $m<2$,
the proposition now follows by taking $m\nearrow2$ and using
(\ref{eq5.10}), (\ref{eq5.6}), and the inferior part of
(\ref{eq4.13}). We note that we obtain $F^+_{\lambda_2}(T)\le
const\cdot e^{PT}$ with $P$ as in (\ref{eq5.10}).\qed
\bigskip

It immediately follows that the other sectional curvature,
$\lambda_1$, is bounded below:

%PROPOSITION 5.3%%%%%%%%%%%%%%%%%%%%%%%%%%%%%%%%%%%%%%%%%
\bigskip
\noindent {\bf Corollary 5.3.} {\it Assume that (\ref{eq5.1}) holds
for all $T>0$. Then there is a function
$F^-_{\lambda_1}:[0,\infty)\times [0,\infty)$ such that}
\begin{equation} \lambda_1(t,r)\ge
-F^-_{\lambda_1}(T)\ {\rm whenever\ } 0\le t \le T\ .\label{eq5.12}
\end{equation}

\bigskip
\noindent{\bf Proof.} In rotational symmetry, we have
\begin{equation}
R=2(n-1)\lambda_1+(n-1)(n-2)\lambda_2\ . \label{eq5.13}
\end{equation}
The result then follows from (\ref{eq4.9}) and Proposition 5.2, and
indeed $F^-_{\lambda_1}(T)\le const\cdot e^{PT}$ with $P$ as in
(\ref{eq5.10}).\qed
\bigskip

\subsection{Bounding the Hessian of $u$}

Finally we seek an upper bound for $\lambda_1$. To find it, we must
first bound the second $r$-derivative of $u$.

%PROPOSITION 5.4%%%%%%%%%%%%%%%%%%%%%%%%%%%%%%%%%%%%%%%%%
\bigskip
\noindent {\bf Proposition 5.4.}  {\sl Assume that (\ref{eq5.1})
holds for all $T>0$. Then there is a function
$F^+_{|z'|}:[0,\infty)\times [0,\infty)$ such that}
\begin{equation}
\left \vert z'(t,r) \right \vert \le F^+_{\left \vert z' \right \vert }(T)
\ {\rm whenever\ } 0\le t \le T\  . \label{eq5.14}
\end{equation}

\bigskip
\noindent{\bf Proof.}  Recall that in rotational symmetry we have
$|z|=\left \vert \frac1f \frac{\partial u}{\partial r}\right \vert
=|\nabla u|$, with $|z'|=\left \vert \frac{\partial}{\partial
r}\left(\frac1f \frac{\partial u}{\partial r}\right)\right \vert
=\frac{\partial}{\partial r}|\nabla u|$.  The evolution equation for
$(z')^2$ can be derived from (\ref{eq3.11}, \ref{eq3.12}) and is
given by
\begin{eqnarray}
\frac{\partial}{\partial t}\left [ (z')^2 \right ]
%&=&\frac{1}{f^2}\frac{\partial^2(z')^2}{\partial r^2}
%+\left( \frac{1}{rf^2}+\frac{n-2}{r}-\frac{2}{f^3}\frac{\partial f}{\partial r}\right)
%\frac{\partial (z')^2}{\partial r}
%-\frac{2}{f^2}\left(\frac{\partial z'}{\partial r}\right)^2\nonumber \\
%&&-2\left (\frac{R}{n-1}+\frac{2(n-1)}{r^2f^2}+3k_n^2z^2\right )
%(z')^2 \nonumber\\
%&&-2\left [R+\left(\frac{n-1}{r^2}\right)\left(\frac{n}{f^2}-(n-2)\right)\right ]\frac{z}{r}z'\nonumber\\
&=&\frac{1}{f^2}\frac{\partial^2}{\partial r^2}\left [ (z')^2 \right ]
+\left( \frac{1}{rf^2}+\frac{n-2}{r}-\frac{2}{f^3}
\frac{\partial f}{\partial r}\right) \frac{\partial}{\partial r}
\left [ (z')^2 \right ]
-\frac{2}{f^2}\left(\frac{\partial z'}{\partial r}\right)^2\nonumber\\
&&-2\left(R+\frac{n(n-1)}{r^2f^2}\right)\left(\frac{(z')^2}{n-1}
+\frac{z}{r}z'\right)-6k_n^2z^2(z')^2\nonumber\\
&&-\frac{2}{r^2}(n-1)(n-2)\left(\frac{(z')^2}{(n-1)f^2}
-\frac{z}{r}z'\right) \label{eq5.15}
\end{eqnarray}

Set $0\le t \le T$ and define
\begin{equation}
{\cal Z}:=e^{-2\kappa t}(z')^2\ , \label{eq5.16}
\end{equation}
where the constant $\kappa >0$ will be chosen below. Then from
(\ref{eq5.15}) we compute that
\begin{eqnarray}
\frac{\partial {\cal Z}}{\partial t}&=&\frac{1}{f^2}
\frac{\partial^2 {\cal Z}}{\partial r^2} +\left ( \frac{1}{rf^2}
+\frac{n-2}{r}-\frac{2}{f^3}\frac{\partial f}{\partial r} \right )
\frac{\partial {\cal Z}}{\partial r}-\frac{1}{2f^2{\cal Z}}\left (
\frac{\partial {\cal Z}}{\partial r} \right )^2\nonumber \\
&& -2\left ( R +\frac{n(n-1)}{r^2f^2} \right ) \left [ \frac{\cal Z}{n-1}
+\frac{z}{r} e^{-\kappa t}\sqrt{\cal Z} \right ]\label{eq5.17}\\
&&-\frac{2}{r^2}(n-1)(n-2)\left [ \frac{\cal Z}{(n-1)f^2}-\frac{z}{r}
e^{-\kappa t}\sqrt{\cal Z} \right ]-2\left [ 3k^2_n z^2+\kappa\right ]
{\cal Z} \ . \nonumber
\end{eqnarray}

By asymptotic flatness, ${\cal Z}\to 0$ as $r\to\infty$ for an
$t<T$. Since ${\cal Z}\ge 0$ by definition, either ${\cal Z}=z'=0$
everywhere or ${\cal Z}$ has a positive maximum. Say the maximum of
${\cal Z}$ occurs at a spacetime point $q=(t,r)$.

\begin{enumerate}
\item[{\it (i)}] {\it Case $R=:R_q\ge 0$ at $q$ and $r(q)\neq
    0$}: From (\ref{eq5.17}) it follows that at $q$ we must have
\begin{eqnarray}
0\le\frac{\partial {\cal Z}}{\partial t}&\le& -2 \left [ R_q
+\frac{n(n-1)}{r^2f^2} \right ] \left [ \frac{\cal Z}{n-1}+\frac{z}{r}
e^{-\kappa t}\sqrt{\cal Z} \right ]\nonumber\\
&&-\frac{2}{r^2}(n-1)(n-2)\left [ \frac{\cal Z}{(n-1)f^2}-\frac{z}{r}
e^{-\kappa t}\sqrt{\cal Z} \right ]\ , \label{eq5.18}
\end{eqnarray}
and so we must have that
\begin{equation}
\sqrt{\cal Z}\le -(n-1)\frac{z}{r}e^{-\kappa t}
\le (n-1)F^+_{\zeta}(T)e^{-\kappa t} \label{eq5.19}
\end{equation}
if $n=2$, and if $n>2$ then either (\ref{eq5.19}) must hold or
\begin{equation}
\sqrt{\cal Z}\le (n-1)f^2\frac{z}{r}e^{-\kappa t}
\le F^+_{\zeta}(T)(C^+_f)^2(1+t)^{2p}e^{-\kappa t} \label{eq5.20}
\end{equation}
must hold instead, using (\ref{eq5.1}) and (\ref{eq4.13}) of
Proposition 4.5 (where $p$ is defined) and (\ref{eq4.5}). In
either case, ${\cal Z}$ and, thus, $z'$ are bounded above at any
$t\ge 0$.

\item[{\it (ii)}] {\it Case $R< 0$ at $q$ and $r(q)\neq 0$}: We
    further assume that neither (\ref{eq5.19}) nor
    (\ref{eq5.20}) holds at the maximum, since otherwise we
    would have an upper bound on ${\cal Z}$. Then the terms in
    square brackets in (\ref{eq5.17}) are non-negative so where
    they are multiplied by negative coefficients in
    (\ref{eq5.17}) we can drop them and obtain that
\begin{equation}
0\le\frac{\partial {\cal Z}}{\partial t} \le -2R_q\left ( \frac{\cal Z}{n-1}
+\frac{z}{r}e^{-\kappa t}\sqrt{\cal Z} \right ) - 2\kappa{\cal Z}
\label{eq5.21}
\end{equation}
at the point $q$ where ${\cal Z}$ takes its maximum. Choose
$\kappa$ such that $\kappa+\frac{R_q}{n-1}>0$. For example,
choose
\begin{equation}
\kappa=1+\frac{\left \vert C_S^-\right \vert}{n-1}\ , \label{eq5.22}
\end{equation}
where $C_S^-$ is a lower bound for $R$ (cf (\ref{eq4.9})). Then
(\ref{eq5.21}) yields
\begin{equation}
\sqrt{\cal Z}\le -R_q\frac{z}{r}e^{-\kappa t} \le
|C_s^-|F^+_{\zeta}(T)e^{-\kappa t}\ . \label{eq5.23}
\end{equation}

\item[{\it (iii)}] {\it Maximum occurs at $r=r(q)=0$}: By local
    existence, $f$ and $z$ and their spatial derivatives are
    bounded for $0\le t \le T$, for any $T<\tau=$ maximal time
    of existence. Thus the same is true for ${\cal Z}$ and for
    its first time derivative. Examining the behaviour of
    coefficients in (\ref{eq5.17}) as $r\to 0$, keeping in mind
    that $z/r$ is bounded, we see that this implies that
\begin{equation}
\frac{n}{f^2}\left [ \frac{\cal Z}{n-1}+\frac{z}{r}
e^{-\kappa t}\sqrt{\cal Z} \right ]+(n-2)\left [ \frac{\cal Z}{(n-1)f^2}
-\frac{z}{r} e^{-\kappa t}\sqrt{\cal Z} \right ] \in {\cal O}(r^2)
\label{eq5.24}
\end{equation}
as $r\to 0$ for all $0\le t\le T$. Taking the limit as $r\to 0$
of (\ref{eq5.24}), we obtain either
\begin{equation}
\sqrt{{\cal Z}(t,0)}=0\ , \label{eq5.25}
\end{equation}
or
\begin{equation}
\sqrt{{\cal Z}(t,0)}=\frac12 f^2
\left ( n-2-\frac{n}{f^2} \right ) \frac{z}{r}e^{-\kappa t} \ .
\label{eq5.26}
\end{equation}

\item[{\it (iv)}] {\it Maximum occurs at $t=0$}: Then from
    (\ref{eq5.16}) we get
\begin{equation}
\sqrt{{\cal Z}(0,r)}\equiv \left \vert z'(0,t) \right \vert\le const\ ,
\label{eq5.27}
\end{equation}
since the initial data for $z'$ is bounded.

\end{enumerate}

At least one of the bounds given by (\ref{eq5.19}), (\ref{eq5.20}),
(\ref{eq5.23}), (\ref{eq5.25}--\ref{eq5.27}) must hold and so, using
the definition (\ref{eq5.16}) with $\kappa$ given by (\ref{eq5.22}),
we obtain (\ref{eq5.14}). \qed

\bigskip

\subsection{An upper bound on transverse sectional curvature}

\noindent With the Hessian bound in hand, the following estimate
then gives the desired upper bound for $\lambda_1$.

%PROPOSITION 5.5%%%%%%%%%%%%%%%%%%%%%%%%%%%%%%%%%%%%%%%%%
\bigskip
\noindent {\bf Proposition 5.5.} {\sl Assume that (\ref{eq5.1})
holds for all $T>0$. Then there is a function
$F^+_{\lambda_1}:[0,\infty)\times [0,\infty)$ such that}
\begin{equation}
\lambda_1(t,r)\le F^+_{\lambda_1}(T) \ {\rm whenever\ } 0\le t \le T
\ . \label{eq5.28}
\end{equation}
\bigskip

\noindent{\bf Proof.}
We first define
\begin{equation}
y(t,r):= \begin{cases} \frac12 r\frac{\partial}{\partial r}
\left[\frac{1}{r^2}\left(\frac{1}{f}-1\right)\right]
& \text{\ for\ } r>0\ ,\\
0& \text{\ for\ } r=0 \ , \end{cases}\label{eq5.29}
\end{equation}
and we note that
\begin{equation}
y= \frac{f}{(1+f)}\lambda_2 - \frac12 f\lambda_1 \text{\ for\ }
r>0\ , \label{eq5.30}
\end{equation}
so we seek a lower bound for $y$. To see that $y$ is continuous at
$r=0$, use the Bianchi identity (\ref{eq3.16}) to write
(\ref{eq5.30}) as
\begin{equation}
y= -\frac{rf}{2(1+f)}\frac{\partial \lambda_2}{\partial r}
+ \frac{f(1-f)}{2(1+f)}\lambda_1\ . \label{eq5.31}
\end{equation}
Since $1-f\in {\cal O}(r^2)$ and $\frac{\partial \lambda_2}{\partial
r} \in {\cal O}(r)$ as $r\to 0$, we see that $y\to 0$ as $r\to 0$.

Computing from (\ref{eq3.11}) and (\ref{eq5.29}), we see that for
$r>0$ $y$ obeys
\begin{eqnarray}
\frac{\partial y}{\partial t}
&=&\frac{1}{f^2}\frac{\partial^2y}{\partial r^2}
+\left(\frac{\alpha}{r}\right)\frac{\partial y}{\partial r}
+\frac{8}{f}y^2+\frac{1}{r^2}\left [
\beta y+\gamma \right ]-k_n^2z^2y+k_n^2\frac{z^2}{r^2}
-\frac{k_n^2}{f}\frac{z}{r} z' \ , \label{eq5.32}\\
\alpha&:=&
\frac{4r^2}{f}y+\frac{5}{f^2}-\frac{4}{f}+n-2 \ ,\label{eq5.33} \\
\beta&:=&\frac{4}{f^2} -\frac{8}{f} + (n-2)\left ( 1-\frac{3}{f^2} \right )
\ , \label{eq5.34} \\
\gamma&:=&\frac{(n-2)}{r^2} \left ( 1-\frac{1}{f}\right )^3
\ . \label{eq5.35}
\end{eqnarray}
Using the definition of $\lambda_2$, we simplify $\beta$ and
$\gamma$ as follows:
\begin{eqnarray}
\frac{\beta}{r^2}&=&\left ( n-2-\frac{4}{(1+f)} \right )
\lambda_2-\frac{4}{r^2f} -\frac{2(n-2)}{r^2f^2}\nonumber\\
&\le& \left ( n-2-\frac{4}{(1+f)} \right ) \lambda_2\ , \label{eq5.36}\\
\frac{\gamma}{r^2}&=& (n-2) \frac{f(f-1)}{(f+1)^2}\lambda_2^2
\ge -\frac{(n-2)f}{(f+1)^2}\lambda_2^2\ . \label{eq5.37}
\end{eqnarray}
Then whenever $y<0$ and $r>0$ we have that (\ref{eq5.32}) yields
\begin{eqnarray}
\frac{\partial y}{\partial t}
&\ge&\frac{1}{f^2}\frac{\partial^2y}{\partial r^2}
+\left(\frac{\alpha}{r}\right)\frac{\partial y}{\partial r}
+\frac{8}{f}y^2
+\left ( n-2-\frac{4}{(1+f)} \right ) \lambda_2 y\nonumber \\
&&- \frac{(n-2)f}{(1+f)^2}\lambda_2^2
-\frac{k_n^2}{f}\frac{z}{r} z' \ . \label{eq5.38}
\end{eqnarray}
In particular, we work as usual on the domain $A_{\epsilon}(T):=
[0,T]\times \left [ \epsilon, \frac{1}{\epsilon}\right ]$ and then
observe immediately from (\ref{eq5.38}) that, if $y$ takes its
minimum in the parabolic interior, then $y(t,r)$ is bounded below by
\begin{eqnarray}
y(t,r)&\ge&-\frac12 \left ( n-2-\frac{4}{(1+f)} \right )\frac{f\lambda_2}{8}
\nonumber \\
&&- \left [ \left ( n-2-\frac{4}{(1+f)} \right )^2\frac{f^2\lambda_2^2}{64}
+ \frac{(n-2)f^2}{8(1+f)^2}\lambda_2^2+\frac{k_n^2}{8}\frac{z}{r} z'
\right ]^{1/2} \ , \label{eq5.39}
\end{eqnarray}
and then the bounds (\ref{eq4.13}), (\ref{eq4.37}), (\ref{eq5.1}),
(\ref{eq5.3}), and (\ref{eq5.14}) on the quantities appearing on the
right-hand side prove the proposition.

On the other hand, $y$ could have its minimum on the parabolic
boundary of $A_{\epsilon}(T)$. If the minimum occurs at $t=0$, then
$y$ is bounded below by $\min \left \{ 0, \inf_r \{y(0,r)\} \right
\}$, again proving the proposition. If, however, the minimum occurs
at $r=\epsilon$ or $r=1/\epsilon$, then if we choose $\epsilon$
small enough this minimum would approach zero (since
$y(t,\epsilon)\to 0$ as $r\to 0$ by the argument at the start of the
proof, and $y(t,1/\epsilon)\to 0$ by asymptotic flatness, as seen
from, say, (\ref{eq5.30})). \qed

\section{Proofs of Theorems 1.3 and 1.5}
\setcounter{equation}{0}

\bigskip
\noindent{\bf Proof of 1.3.} By assumption, there is no minimal
hypersphere at $t=0$. By (\ref{eq4.10}) and (\ref{eq4.13}), then no
minimal hypersphere can form at any $t\in [0,T_M)$.

Furthermore, $\lambda_2$ is bounded below for all $t>0$ (Proposition
4.8). Now assume that equation (\ref{eq5.1}) holds. Then for all
$t>0$, $\lambda_2$ is bounded above (Proposition 5.1) and
$\lambda_1$ is bounded below (Corollary 5.3) and above (Proposition
5.5). Thus the maximal time of existence is $T_M=\infty$
(Proposition 3.4).

In particular, if $n=2$, then (\ref{eq5.1}) holds (Proposition 4.9).
\qed

\bigskip
\noindent{\bf Proof of 1.5.} There are no closed geodesics of
$g(t)$, for if there were then by rotational symmetry such a
geodesic would necessarily lie on a minimal hypersphere, and we have
shown that there are none of those. Then by standard results
(\cite{Petersen}, paragraph 6.6.1), the injectivity radius of the
manifold at time $t$ will be equal to the conjugate radius and thus
bounded below by $\pi/\sqrt{\sup_r |{\rm Riem}(t,r)|}$. These facts
are diffeomorphism invariant and so apply equally to ${\bar g}(t)$
(see (\ref{eq3.2}, \ref{eq3.5})). Choosing an essential blow-up
sequence for ${\bar g}, {\bar u})$ and rescaling as in
(\ref{eq1.15}), then along each rescaled flow $(g_{(k)},u_{(k)})$
the injectivity radius of $g_{(k)}(s)$ is uniformly (in $s$ and in
$k$) bounded below by $\pi/\sqrt{C}$.

In view of equation (\ref{eq1.9}), define
$X^i_{(k)}:=-g_{(k)}^{ij}\nabla^{(k)} u_{(k)}$. Note that
\begin{eqnarray}
\left \vert X_{(k)}\right \vert^2&:=&\left [ g_{(k)}^{ij}\nabla_i^{(k)} u_{(k)}
\nabla_j^{(k)} u_{(k)}
\right ]_s = \left [ \frac{1}{B_k}{\bar g}^{ij}{\bar \nabla}_i {\bar u}
{\bar \nabla_j} {\bar u}
\right ]_{(t_k+s/B_k)}\nonumber \\
&\le& \frac{const}{B_k(1+t_k+s/B_k)}\nonumber \\
&\to& 0 \label{eq6.1}
\end{eqnarray}
In particular, for each $k$ the diffeomorphisms generated by
$X_{(k)}$ are defined (and, indeed, getting smaller). Thus, we can
use the correspondence between List's flow and Ricci flow
(\ref{eq1.7}--\ref{eq1.9}) to express the sequence $\left (
g_{(k)}(s),u_{(k)}(s) \right )$ as a sequence of Ricci flows
$G_{(k)}(s)$ in $(n+1)$-dimensions. Because
$\frac{\partial}{\partial \tau}$ (cf (\ref{eq1.8})) is a Killing
vector field, the injectivity radius remains bounded below by
$\pi/\sqrt{C}$. Let $x_k:=(\tau_k,r_k)$ and choose $\tau_k=0$ (since
$\frac{\partial}{\partial \tau}$ is a Killing vector field, the
choice is irrelevant). By a theorem of Hamilton \cite{Hamilton}, the
pointed sequence $(M,G_{(k)}(s),x_k)$ converges to a complete
pointed Ricci flow $(M,G(s),x)$. The domain of $s$ is the limit of
the intervals $[-B_k(1+t_k),0]$ and is thus $(-\infty, 0]$, so $G$
is an ancient solution of Ricci flow. The injectivity radius of
$(M,G(s))$ is bounded below (uniformly in $s$) at $x$.

We see from (\ref{eq6.1}) that $u_{(\infty)}$ is constant in $r$. It
is constant in $\tau=x^0$ by assumption and then is constant in the
flow time $s$ (equivalently, in $t$) as well by the asymptotic
condition $u(t,r)\to const$ as $r\to\infty$. It follows that the
Ricci flow for the limit metric $G(s)$ in $(n+1)$-dimensions is
trivial in the $\tau$ direction, and splits as an ancient Ricci flow
for $g$ (the induced metric for $\tau=0$) in $n$-dimensions,
together with the equation $u=const$.

Since $\lambda_2(t)\ge -const$ for the flow (Prop 4.8) of $(g,u)$
(and thus for the unrescaled flow of $({\bar g}, {\bar u})$ since
the condition is natural with respect to diffeomorphisms) and since
rescaling divides $\lambda_2$ by the maximum of the norm of the
curvature, the limit of rescaled flows is a flow with
$\lambda_2(s)\ge 0$. That is, the limit flow has nonnegative
sectional curvature in planes tangent to the orbits of symmetry. By
Theorem 2.4 of \cite{Chen}, any ancient, complete, 3-dimensional
solution of Ricci flow has nonnegative sectional curvatures in all
planes, thus including radial planes as well as tangential planes
when $n=3$. Chen also observed that any ancient, complete flow has
$R\ge 0$, but here we can see this directly by the same argument as
with $\lambda_2$, since $R\ge -const$ along the original flow
(Proposition 4.4). \qed

\bigskip

\noindent{\bf Acknowledgments.} This work was carried out during
visits of TAO to the University of Alberta, EW to Monash University,
and all three authors to the Banff International Research Station.
TAO and EW also thank the Arnold Sommerfeld Centre of the LMU Munich
and the organizers of the Workshop on Field Theory and Geometric
Flows. EW thanks the Ennio di Giorgi Centre, Scuola Normale
Superiore di Pisa, and the organizers of the Pisa workshop on
Geometric Flows in Mathematics and Theoretical Physics. We thank
Maria Athanassenas and Reto M\"uller for discussions. The work was
partially supported by a Discovery Grant to EW from the Natural
Sciences and Engineering Research Council of Canada.

\end{document}